\documentclass[11pt]{amsart}
\usepackage{amssymb,mathtools,fullpage,aliascnt,eucal}
\usepackage[dvipsnames]{xcolor}
\usepackage[hypertexnames=false]{hyperref}
\usepackage[capitalize,nameinlink,noabbrev,nosort]{cleveref}
\hypersetup{
  colorlinks=true,
  allcolors=RedViolet
}
\DeclareFontEncoding{OT2}{}{}
\DeclareFontSubstitution{OT2}{cmr}{m}{n}
\DeclareFontFamily{OT2}{cmr}{\hyphenchar\font45 }
\DeclareFontShape{OT2}{cmr}{m}{n}{
<5><6><7><8><9>gen*wncyr
<10><10.95><12><14.4><17.28><20.74><24.88>wncyr10}{}
\DeclareFontShape{OT2}{cmr}{b}{n}{
<5><6><7><8><9>gen*wncyb
<10><10.95><12><14.4><17.28><20.74><24.88>wncyb10}{}
\DeclareMathAlphabet{\mathcyr}{OT2}{cmr}{m}{n}
\SetMathAlphabet{\mathcyr}{bold}{OT2}{cmr}{b}{n}
\newcounter{nameditem}

\newcommand{\nameditem}[1]{%
  \item[(#1)]%
  \renewcommand{\thenameditem}{(#1)}%
  \refstepcounter{nameditem}%
}
\theoremstyle{plain}
\newtheorem{theorem}{Theorem}[section]

\newtheorem{maintheorem}{Theorem}

\crefname{maintheorem}{Theorem}{Theorems}
\Crefname{maintheorem}{Theorem}{Theorems}

\newaliascnt{proposition}{theorem}
\newtheorem{proposition}[proposition]{Proposition}
\aliascntresetthe{proposition}

\newaliascnt{conjecture}{theorem}
\newtheorem{conjecture}[conjecture]{Conjecture}
\aliascntresetthe{conjecture}
\crefname{conjecture}{Conjecture}{Conjectures}
\Crefname{conjecture}{Conjecture}{Conjectures}

\newaliascnt{lemma}{theorem}
\newtheorem{lemma}[lemma]{Lemma}
\aliascntresetthe{lemma}

\newaliascnt{corollary}{theorem}
\newtheorem{corollary}[corollary]{Corollary}
\aliascntresetthe{corollary}

\theoremstyle{definition}
\newaliascnt{definition}{theorem}
\newtheorem{definition}[definition]{Definition}
\aliascntresetthe{definition}

\theoremstyle{remark}
\newaliascnt{remark}{theorem}
\newtheorem{remark}[remark]{Remark}
\aliascntresetthe{remark}
\numberwithin{equation}{section}
\newcommand{\calA}{\mathcal{A}}
\newcommand{\rbl}{\mathrm{bl}}
\newcommand{\block}{\mathrm{block}}
\newcommand{\bbC}{\mathbb{C}}
\newcommand{\rcf}{\mathrm{cf}}
\newcommand{\rCF}{\mathrm{CF}}
\newcommand{\calD}{\mathcal{D}}
\newcommand{\dep}{\mathrm{dep}}
\newcommand{\rDe}{\mathrm{De}}
\newcommand{\bbF}{\mathbb{F}}
\DeclareMathOperator{\GRT}{GRT}
\newcommand{\calH}{\mathcal{H}}
\newcommand{\rH}{\mathrm{H}}
\newcommand{\bbI}{\mathbb{I}}
\DeclareMathOperator{\id}{id}
\newcommand{\calI}{\mathcal{I}}
\newcommand{\bk}{\boldsymbol{k}}
\newcommand{\rKZ}{\mathrm{KZ}}
\newcommand{\bl}{\boldsymbol{l}}
\newcommand{\rL}{\mathrm{L}}
\newcommand{\bbQ}{\mathbb{Q}}
\newcommand{\bbR}{\mathbb{R}}
\newcommand{\calR}{\mathcal{R}}
\newcommand{\rR}{\mathrm{R}}
\newcommand{\reg}{\mathrm{reg}}
\newcommand{\calS}{\mathcal{S}}
\newcommand{\sh}{\mathbin{\mathcyr{sh}}}
\DeclareMathOperator{\Span}{span}
\newcommand{\sfw}{\mathsf{w}}
\newcommand{\bbZ}{\mathbb{Z}}
\newcommand{\calZ}{\mathcal{Z}}
\newcommand{\rZ}{\mathrm{Z}}
\newcommand{\sfz}{\mathsf{z}}
\newcommand{\fps}[1]{\langle\!\langle #1\rangle\!\rangle}
\newcommand{\bh}{\boldsymbol{h}}%
\begin{document}
\title{On Zhao's conjectural basis for finite multiple zeta values}
\author{Shin-ichiro Seki}
\thanks{This research was supported by JSPS KAKENHI Grant Number JP26K06734.}
\address{Nagahama Institute of Bio-Science and Technology, 1266, Tamura, Nagahama, Shiga, 526-0829, Japan}
\email{s\_seki@nagahama-i-bio.ac.jp}
\begin{abstract}
We prove a conjecture by J. Zhao which states that every finite multiple zeta value associated with a nonempty index is a $\bbQ$-linear combination of $\zeta^{}_{\calA}(k_1,\dots, k_r, 2, 1)$, where $k_i\in \{2, 3\}$.
Our proof proceeds by constructing a basis for the space of formal multiple zeta values defined by the confluence relations.
This approach also reduces the Kaneko--Zagier conjecture (apart from compatibility with the depth filtration) to linear independence conjectures for ordinary and finite multiple zeta values, without using motivic theory.
\end{abstract}
\maketitle
\section{Introduction}
\subsection{Main results}
An \emph{index} is a finite, possibly empty, sequence of positive integers, and we denote the set of all indices by $\bbI$.
The empty index is denoted by $\varnothing$.
A nonempty index $\bk=(k_1,\dots,k_r)$ is called \emph{admissible} if $k_r\geq2$, and the empty index is also regarded as admissible.
The \emph{weight} and \emph{depth} of $\bk$ are defined by $\lvert\bk\rvert\coloneqq k_1+\cdots+k_r$ and $\dep(\bk)\coloneqq r$, respectively.
We set $\bbI_k\coloneqq \{\bk\in\bbI : \lvert\bk\rvert=k\}$ and $\bbI_{k,r}\coloneqq \{\bk\in\bbI_k : \dep(\bk)=r\}$.
It is admissible, and we set $\lvert\varnothing\rvert=\dep(\varnothing)=0$.
We set
\begin{align*}
\bbI_k^{\rH}&\coloneqq \{\bk\in\bbI_k : \bk=(k_1,\dots,k_r),\;r\in\bbZ_{\geq 0},\; k_1,\dots,k_r\in\{2,3\}\},\\
\bbI_k^{\rZ}&\coloneqq \{\bk\in\bbI_k : \bk=(k_1,\dots,k_r,2,1),\;
r\in\bbZ_{\geq 0},\; k_1,\dots,k_r\in\{2,3\}\}.
\end{align*}
The elements of $\bbI_k^{\rH}$ and $\bbI_k^{\rZ}$ are called \emph{Hoffman indices} and \emph{Zhao indices} of weight $k$, respectively.
Zhao uses the opposite convention for the order of summation, so his proposed indices appear in reversed order in our notation.
Let $(d_k)_{k\in\bbZ_{\geq 0}}$ be the sequence defined by the generating function
\[
\sum_{k=0}^{\infty}d_kX^k=\frac{1}{1-X^2-X^3}.
\]
Then $\#\bbI_k^{\rH}=d_k$ for all $k\geq 0$, and $\#\bbI_k^{\rZ}=d_{k-3}$ for all $k\geq 3$.

For an admissible index $\bk=(k_1,\dots,k_r)$, the \emph{multiple zeta value} $\zeta(\bk)\in\bbR$ is defined by
\[
\zeta(k_1,\dots,k_r)\coloneqq \sum_{0<n_1<\cdots<n_r}\frac{1}{n_1^{k_1}\cdots n_r^{k_r}}.
\]
We set $\zeta(\varnothing)=1$.
Similar conventions will be used without further mention.
Let $\calZ_k$ denote the $\bbQ$-vector space spanned by $\{\zeta(\bk) : \bk\in\bbI_k \text{ is admissible}\}$.
Zagier \cite{Zagier1994} conjectured that $\dim_{\bbQ}\calZ_k=d_k$.
The upper bound in this conjecture has been established as follows.
\begin{theorem}[Terasoma \cite{Terasoma2002}, Deligne--Goncharov \cite{DeligneGoncharov2005}]\label{thm:TerasomaDeligneGoncharov}
For all $k\geq 0$, $\dim_{\bbQ}\calZ_k\leq d_k$.
\end{theorem}
More specifically, Hoffman \cite{Hoffman1997} proposed a conjectural basis consisting of multiple zeta values with Hoffman indices.
Brown subsequently proved that these values generate each space $\calZ_k$.
\begin{theorem}[Brown \cite{Brown2012}]\label{thm:Brown}
Every multiple zeta value is a $\bbQ$-linear combination of multiple zeta values with Hoffman indices of the same weight.
\end{theorem}
In this paper, we prove an analogue of Brown's theorem for finite multiple zeta values.
Consider the ring
\[
\calA\coloneqq\Biggl(\prod_{p:\text{ prime}}\bbZ/p\bbZ\Biggr) \Bigg/ \Biggl(\bigoplus_{p:\text{ prime}}\bbZ/p\bbZ\Biggr).
\]
The ring $\calA$ is naturally a $\bbQ$-algebra: a rational number is reduced modulo every prime not dividing its denominator, and the remaining components may be chosen arbitrarily.
For an index $\bk=(k_1,\dots,k_r)$, the \emph{finite multiple zeta value} $\zeta^{}_{\calA}(\bk)\in\calA$ is defined by $\zeta^{}_{\calA}(\bk)\coloneqq (\zeta^{}_{<p}(\bk)\bmod p)_p\in\calA$, where
\[
\zeta^{}_{<p}(k_1,\dots,k_r)\coloneqq\sum_{0<n_1<\cdots<n_r<p}\frac{1}{n_1^{k_1}\cdots n_r^{k_r}}.
\]
Let $\calZ_{\calA,k}$ denote the $\bbQ$-vector space spanned by $\{\zeta^{}_{\calA}(\bk) : \bk\in\bbI_k\}$.
Zagier and Zhao independently conjectured that $\dim_{\bbQ}\calZ_{\calA,k}=d_{k-3}$.
We set $d_{-1}=d_{-2}=0$ and $d_{-3}=1$.
The following counterpart of \cref{thm:TerasomaDeligneGoncharov} is known.
\begin{theorem}[Akagi--Hirose--Yasuda \cite{AkagiHiroseYasuda}]
For all $k\geq 0$, $\dim_{\bbQ}\calZ_{\calA, k}\leq d_{k-3}$. 
\end{theorem}
Akagi, Hirose, and Yasuda conjectured an identity between Deligne's $p$-adic multiple zeta values and multiple harmonic sums, which was subsequently proved by Jarossay \cite{Jarossay}.
Using this identity, they announced a proof of this upper bound, but their paper \cite{AkagiHiroseYasuda} has not been made publicly available as of this writing.

More precisely, Zhao conjectured the following:
\begin{conjecture}[Zhao {\cite[Conjecture~8.1 (ii)]{Zhao}}]
The $\bbQ$-vector space $\calZ_{\calA,k}$ has a basis consisting of the elements of the set $\{\zeta^{}_{\calA}(\bk) : \bk\in\bbI^{\rZ}_k\}$ for all $k\geq 3$.
\end{conjecture}
In this paper, we prove that the family proposed in the above conjecture generates the space.
\begin{maintheorem}\label{thm:mainA}
Every finite multiple zeta value for a nonempty index is a $\bbQ$-linear combination of finite multiple zeta values with Zhao indices of the same weight.
\end{maintheorem}
Unlike ordinary multiple zeta values, which take values in the real numbers, finite multiple zeta values take values in the ring $\calA$, which simultaneously encodes reductions modulo $p$ for all sufficiently large primes $p$.
Despite this substantial difference in their target spaces, the spaces of multiple zeta values and finite multiple zeta values exhibit many similarities.
For instance, their dimension conjectures are formulated in terms of the same sequence $(d_k)$, and their conjectural bases have similar structures.
However, as in the case of ordinary multiple zeta values, the lower bound predicted by the dimension conjecture remains unproved.
Such similarities also appear in the explicit relations satisfied by multiple zeta values and their finite counterparts.
Although these relations often share common features, they do not necessarily take exactly the same form.
For example, as an analogue of the Ohno relation for multiple zeta values proved by Ohno \cite{Ohno1999}, Oyama \cite{Oyama2018} established an Ohno-type relation for finite multiple zeta values.
While these similarities are interesting in their own right, the significance of finite multiple zeta values goes beyond their resemblance to ordinary multiple zeta values.
Kaneko and Zagier conjectured that the space of finite multiple zeta values $\calZ_{\calA}\coloneqq\sum_{k\geq 0}\calZ_{\calA,k}$ and the space of multiple zeta values $\calZ\coloneqq\sum_{k\geq 0}\calZ_k$ are explicitly related through the quotient of the latter by the ideal generated by $\zeta(2)$.
The stuffle product formulas show that $\calZ$ and $\calZ_{\calA}$ are $\bbQ$-algebras under their usual multiplication.

To state this conjecture, we introduce symmetric multiple zeta values together with the Hoffman algebra setup that will also be needed later.
Let $\calH\coloneqq\bbQ\langle e_0,e_1\rangle$ be the noncommutative polynomial ring in the two variables $e_0$ and $e_1$.
We consider the subalgebras $\calH^0\subset\calH^1\subset\calH$ given by $\calH^0\coloneqq\bbQ+e_1\calH e_0$ and $\calH^1\coloneqq\bbQ+e_1\calH$.
For an index $\bk=(k_1,\dots,k_r)\in\bbI$, write $\sfw(\bk)\coloneqq e_1e_0^{k_1-1}\cdots e_1e_0^{k_r-1}$, with $\sfw(\varnothing)=1$.
We denote the usual shuffle product by $\sh$.
Under the standard algebra isomorphism $(\calH,\sh)\simeq(\calH^0,\sh)[e_0,e_1]$, evaluation at $e_0=e_1=0$ on the right-hand side gives the regularization map $\reg_{\sh}\colon\calH\to\calH^0$.
Let $Z\colon\calH^0\to\calZ$ be the $\bbQ$-linear map determined by $Z(\sfw(\bk))=\zeta(\bk)$ for all admissible indices $\bk$.
For any index $\bk\in\bbI$, the shuffle-regularized multiple zeta value $\zeta^{\sh}(\bk)$ is defined by $\zeta^{\sh}(\bk)\coloneqq Z(\reg_{\sh}(\sfw(\bk)))$.
With this notation, for an index $\bk=(k_1,\dots,k_r)\in\bbI$, we set
\[
\zeta^{\sh}_{\calS}(\bk)\coloneqq\sum_{i=0}^r(-1)^{k_{i+1}+\cdots+k_r}\zeta^{\sh}(k_1,\dots,k_i)\zeta^{\sh}(k_r,\dots,k_{i+1})\in\calZ_{\lvert\bk\rvert}.
\]
The \emph{symmetric multiple zeta value} $\zeta^{}_{\calS}(\bk)\in\calZ/\zeta(2)\calZ$ is then defined by $\zeta^{}_{\calS}(\bk)\coloneqq\zeta^{\sh}_{\calS}(\bk)\bmod{\zeta(2)\calZ}$.
\begin{conjecture}[Kaneko--Zagier \cite{KanekoZagier}]\label{conj:KanekoZagier}
There exists a $\bbQ$-algebra isomorphism
\[
\calZ_{\calA}\xrightarrow{\;\sim\;}\calZ/\zeta(2)\calZ
\]
defined by $\zeta^{}_{\calA}(\bk)\mapsto\zeta^{}_{\calS}(\bk)$ for each $\bk\in\bbI$ compatible with the gradings by weight and the filtrations by depth.
Here the depth filtration on $\calZ/\zeta(2)\calZ$ is induced by the usual depth filtration on $\calZ$, while the depth filtration on $\calZ_{\calA}$ is defined so that, for each $d\in\bbZ_{\geq 0}$, its depth-$d$ part is spanned by $1$ and the values $\zeta^{}_{\calA}(\bk)$ with $\dep(\bk)\leq d+1$.
\end{conjecture}
In connection with this conjecture, one may also consider the following finite-dimensional version.
We set $\calZ_m=0$ for $m<0$.
\begin{conjecture}[Finite-dimensional version of the Kaneko--Zagier conjecture]\label{conj:KanekoZagier_weight}
Let $k\geq 2$ be an integer.
There exists a $\bbQ$-linear isomorphism
\[
\calZ_{\calA,k}\xrightarrow{\;\sim\;}\calZ_k/\zeta(2)\calZ_{k-2}
\]
defined by $\zeta^{}_{\calA}(\bk)\mapsto\zeta^{}_{\calS}(\bk)$ for each $\bk\in\bbI_k$.
\end{conjecture}
Note, however, that $\zeta^{}_{\calS}(\bk)$ is used here by an abuse of notation to denote $\zeta^{\sh}_{\calS}(\bk)\bmod{\zeta(2)\calZ_{k-2}}$.
If this conjecture is true, then an analogue of \cref{thm:mainA} for symmetric multiple zeta values would follow from \cref{thm:mainA}.
In fact, we establish this analogue in the present paper.
\begin{maintheorem}\label{thm:mainB}
Every symmetric multiple zeta value for a nonempty index is a $\bbQ$-linear combination of symmetric multiple zeta values with Zhao indices of the same weight.
\end{maintheorem}
In fact, we prove the following stronger theorem, which expresses multiple zeta values in terms of $\zeta(2)$ and $\zeta^{\sh}_{\calS}$-values.
For integers $N\geq r\geq1$ and $\bk=(k_1,\dots,k_r)\in\bbI_{N,r}$, set
\[
\sfz(\bk)\coloneqq(\{2\}^{k_1-1},3,\dots,\{2\}^{k_{r-1}-1},3,\{2\}^{k_r},1)\in\bbI^{\rZ}_{2N+r}.
\]
Every Zhao index arises uniquely in this way.
\begin{maintheorem}\label{thm:mainC}
Let $k$ be a positive integer.
The $\bbQ$-vector space $\calZ_k$ is generated by
the elements of the set
\[
\left\{
  \zeta(2)^a\zeta^{\sh}_{\calS}(\sfz(\bk)) :
  \begin{gathered}
    \bk\in\bbI_{N,r},\quad a,N,r\in\bbZ_{\geq 0},\\
    N\geq r\geq 1,\quad 2a+2N+r=k
  \end{gathered}
\right\}
\cup
\left\{
  \zeta(2)^a :
  a\in\bbZ_{\geq 0},\; 2a=k
\right\}.
\]
\end{maintheorem}
The number of displayed generators, counted with multiplicity, is $d_k$; see \cref{lem:gen=d_k} below.
Thus, this provides a conjectural basis for $\calZ_k$ distinct from the Hoffman basis.

\cref{thm:mainB} follows immediately from \cref{thm:mainC}
by reducing modulo $\zeta(2)$.
\cref{thm:mainC} yields an alternative proof of the following theorem of Yasuda.
For every nonnegative integer $k$, let $\calZ^{\sh}_{\calS,k}$ denote the $\bbQ$-subspace of $\calZ_k$ generated by $\{\zeta^{\sh}_{\calS}(\bk):\bk\in\bbI_k\}$.
\begin{corollary}[Yasuda {\cite[Theorem~6.4]{Yasuda2016}}]\label{thm:Yasuda}
For every nonnegative integer $k$, we have $\calZ^{\sh}_{\calS,k}=\calZ_k$.
\end{corollary}
To deduce $\calZ^{\sh}_{\calS,k}=\calZ_k$ from \cref{thm:mainC}, it suffices to use the inclusion $\zeta(2)\calZ^{\sh}_{\calS,k}\subset\calZ^{\sh}_{\calS,k+2}$; see \cite[Proposition~6.5]{Yasuda2016}.
The resulting proof differs from those of Yasuda \cite{Yasuda2016} and Kina \cite{Kina}.
We do not claim, however, that our proof is more elementary than theirs.
The main point is that our approach establishes the stronger result stated in \cref{thm:mainC}.
The alternative proof of Yasuda's theorem is merely a by-product of this stronger result.

Since both \cref{thm:mainA} and \cref{thm:mainB} follow from \cref{thm:mainF} below, we can construct an explicit family of relations common to finite and symmetric multiple zeta values that suffices to bound the dimensions of both spaces by their conjectured values.
We call the total degree of a word $w\in\calH$ its \emph{length} and denote it by $\lvert w\rvert$.
For a positive integer $k$, let $\calH^1_k$ denote the subspace of $\calH^1$ spanned by the words of length $k$ belonging to $\calH^1$.
Then $\dim_{\bbQ}\calH^1_k=2^{k-1}$.
We define the surjective $\bbQ$-linear maps $Z_{\calA,k}\colon\calH^1_k\twoheadrightarrow\calZ_{\calA,k}$ and $Z_{\calS,k}\colon\calH^1_k\twoheadrightarrow\calZ_k/\zeta(2)\calZ_{k-2}$ by $\sfw(\bk)\mapsto\zeta^{}_{\calA}(\bk)$ and $\sfw(\bk)\mapsto\zeta^{}_{\calS}(\bk)$, respectively.
In this notation, \cref{conj:KanekoZagier_weight} predicts that $\ker Z_{\calA,k}=\ker Z_{\calS,k}$ and, together with the dimension conjectures, that this common kernel has dimension $2^{k-1}-d_{k-3}$.
\begin{maintheorem}\label{thm:mainD}
Let $k$ be a positive integer.
We can explicitly construct a subspace
\[
\calR_k\subset\ker Z_{\calA,k}\cap\ker Z_{\calS,k} \quad \text{with} \quad \dim_{\bbQ}\calR_k=2^{k-1}-d_{k-3}.
\]
Here, the construction is explicit in the sense that there is a deterministic algorithm that, given $k$, computes a basis of $\calR_k$.
See \eqref{eq:def_R_k} below for the definition of $\calR_k$.
\end{maintheorem}
As a corollary of this theorem, we obtain the following.
\begin{maintheorem}\label{thm:mainE}
Assume that the following hypotheses hold:
\begin{description}
\nameditem{Dim $\bbR$}\label{it:dim_MZV} $\dim_{\bbQ}\calZ_k\geq d_k$ for all $k\geq 0$.
\nameditem{Dim $\calA$}\label{it:dim_FMZV} $\dim_{\bbQ}\calZ_{\calA,k}\geq d_{k-3}$ for all $k\geq 3$.
\nameditem{DirSum $\bbR$}\label{it:direct_sum_MZV} (direct sum conjecture for MZVs) $\displaystyle \calZ=\bigoplus_{k\geq 0}\calZ_k$.
\nameditem{DirSum $\calA$}\label{it:direct_sum_FMZV} (direct sum conjecture for FMZVs) $\displaystyle \calZ_{\calA}=\bigoplus_{k\geq 0}\calZ_{\calA,k}$.
\end{description}
Then the Kaneko--Zagier conjecture \cref{conj:KanekoZagier} holds, except for the assertion concerning the depth filtration.
\end{maintheorem}
For the depth filtration, see \cref{rem:depth_KZconj}.
Apart from the explicit construction of $\calR_k$, \cref{thm:mainD} and \cref{thm:mainE} are not new.
Indeed, they can be deduced from the motivic results of Akagi--Hirose--Yasuda \cite{AkagiHiroseYasuda}, which were announced in \cite[Section~2]{Yasuda2018} and are stated more explicitly by Rosen in \cite[Theorem~8.11]{Rosen}.
In the present paper, however, the results are both formulated and proved without using the theory of mixed Tate motives.

The results above arise from a common basis theorem for formal multiple zeta values defined by the confluence relations.
We introduce this formal setting and state the basis theorem in the next subsection.
\subsection{Formal multiple zeta values with respect to the confluence relations}
Let $\calI_{\rCF}$ be the subspace of $\calH^0$ spanned by the confluence relations defined in \cite[Definition~22]{HiroseSato2019}.
Let $L\colon\calH^0\to\calZ$ be the $\bbQ$-linear map determined by $L(\sfw(\bk))=(-1)^{\dep(\bk)}\zeta(\bk)$ for every admissible index $\bk$.
Hirose and Sato proved that multiple zeta values satisfy the confluence relations, that is, $\calI_{\rCF}\subset\ker L$.
See \cite[Theorem~21]{HiroseSato2019}.
It is further conjectured that these relations exhaust all linear relations among multiple zeta values, that is, $\calI_{\rCF}=\ker L$.
See \cite[Conjecture~26]{HiroseSato2019}.

We define the \emph{formal multiple zeta value $\zeta^{\rcf}(\bk)$ with respect to the confluence relations} by
\begin{equation}\label{eq:def_of_zeta^cf}
\zeta^{\rcf}(\bk)\coloneqq(-1)^{\dep(\bk)}\reg_{\sh}(\sfw(\bk))\bmod{\calI_{\rCF}}
\end{equation}
for each index $\bk$, and write $\calZ^{\rcf}\coloneqq\calH^0/\calI_{\rCF}$ for the space spanned by these elements.
For each nonnegative integer $k$, let $\calH^0_k$ denote the $\bbQ$-subspace of $\calH^0$ spanned by the words of length $k$.
Let $\calZ^{\rcf}_k\coloneqq\calH^0_k/(\calI_{\rCF}\cap\calH^0_k)$ denote the subspace of $\calZ^{\rcf}$ generated by the values $\zeta^{\rcf}(\bk)$ with $\bk\in\bbI_k$.
We set $\calZ^{\rcf}_m=0$ for $m<0$.
Since $\calI_{\rCF}$ is a shuffle ideal homogeneous with respect to word length, $\calZ^{\rcf}$ is a commutative $\bbQ$-algebra graded by weight.

Hirose proved the following result for formal multiple zeta values with respect to the confluence relations.
Here, $\bbZ_{(2)}$ denotes the localization of $\bbZ$ at the prime ideal $(2)$; equivalently, it is the subring of $\bbQ$ consisting of rational numbers whose denominators in lowest terms are odd.
\begin{theorem}[Hirose {\cite[Theorem~4]{Hirose}}]\label{thm:Hirose}
Let $k$ be a positive integer.
For every admissible index $\bk\in\bbI_k$, $\zeta^{\rcf}(\bk)$ is a $\bbZ_{(2)}$-linear combination of $\zeta^{\rcf}(\bl)$ with $\bl\in\bbI^{\rH}_k$.
In particular, $\dim_{\bbQ}\calZ^{\rcf}_k\leq d_k$.
\end{theorem}
Since $\calI_{\rCF}\subset\ker L$, the map $L$ induces a surjective $\bbQ$-algebra homomorphism $\calZ^{\rcf}\twoheadrightarrow\calZ$ sending $\zeta^{\rcf}(\bk)$ to $\zeta^{\sh}(\bk)$ for every index $\bk$.
Its restriction to $\calZ^{\rcf}_k$ has image $\calZ_k$.
Applying this evaluation map, we recover \cref{thm:TerasomaDeligneGoncharov} and \cref{thm:Brown}.
Hirose's argument provides more elementary alternative proofs of these theorems without using the theory of motives.
Moreover, his proof gives an explicit algorithm for expressing formal multiple zeta values as $\bbZ_{(2)}$-linear combinations of formal multiple zeta values with Hoffman indices.

In the present paper, we use Hirose's upper bound $\dim_{\bbQ}\calZ^{\rcf}_k\leq d_k$ as a key input, allowing us to avoid motivic arguments.
An alternative route to avoiding motivic methods would be to use the upper bound for the Grothendieck--Teichm\"uller Lie algebra established by OpenAI \cite{OpenAI} prior to Hirose's work, for which Willwacher \cite{Willwacher} subsequently gave a simplified proof.

However, in keeping with the author's preference for an exposition that is as elementary as possible, we not only avoid motivic methods but also make very little use of the language of Lie algebras and Hopf algebras.
From this perspective, basing our argument on Hirose's work is advantageous.
We nevertheless emphasize that, although our approach does not rely on the theory of motives, it draws on several results on associators.
These include Furusho's theorem \cite{Furusho2022}, which establishes the equivalence between the pentagon equation and the confluence relations for commutator group-like formal power series.
In addition, the proof of \cref{thm:mainA} relies on several deep results in $p$-adic theory.

For an index $\bk=(k_1,\dots,k_r)$, we define the formal analogue $\zeta^{\rcf}_{\calS}(\bk)$ of $\zeta^{\sh}_{\calS}(\bk)$ by
\begin{equation}\label{eq:def_of_sym-zeta^cf}
\zeta^{\rcf}_{\calS}(\bk)\coloneqq\sum_{i=0}^r(-1)^{k_{i+1}+\cdots+k_r}\zeta^{\rcf}(k_1,\dots,k_i)\zeta^{\rcf}(k_r,\dots,k_{i+1}).
\end{equation}
The following is the main result of this paper.
\begin{maintheorem}\label{thm:mainF}
Let $k$ be a nonnegative integer.
The $\bbQ$-vector space $\calZ^{\rcf}_k$ has a basis consisting of the elements of the set
\[
\left\{\zeta^{\rcf}(2)^a\zeta^{\rcf}_{\calS}(\sfz(\bk)) :
\begin{gathered}
\bk\in\bbI_{N,r},\quad a,N,r\in\bbZ_{\geq 0},\\
N\geq r\geq 1,\quad 2a+2N+r=k
\end{gathered}
\right\}
\cup
\left\{
\zeta^{\rcf}(2)^a :
a\in\bbZ_{\geq 0},\; 2a=k
\right\}.
\]
\end{maintheorem}
\cref{thm:mainC} follows immediately from \cref{thm:mainF} by applying the evaluation map $\calZ^{\rcf}_k\twoheadrightarrow\calZ_k$.
Note that this evaluation map sends $\zeta^{\rcf}_{\calS}(\bk)$ to $\zeta^{\sh}_{\calS}(\bk)$.

\cref{thm:mainF} yields, in particular, the dimension formula $\dim_{\bbQ}\calZ^{\rcf}_k=d_k$.
This equality could already be deduced by combining Hirose's upper bound with Brown's motivic theorem \cite[Theorem~1.1]{Brown2012}, along the lines of \cite[Section~1.2]{Hirose}.
Our proof, however, does not rely on Brown's theorem or on motivic arguments.
It also follows that the elements of the set $\{\zeta^{\rcf}(\bk) : \bk\in\bbI^{\rH}_k\}$ form a basis of the $\bbQ$-vector space $\calZ^{\rcf}_k$.

The final displayed formula in \cite[Section~1.3]{Hirose} provides an algorithm for expanding $\zeta^{\rcf}(\bk)$ in the Hoffman basis.
Expanding products by the shuffle product and applying this algorithm to the basis elements in \cref{thm:mainF} yields an explicitly computable transition matrix.
Since both families form bases, this matrix is invertible, and its inverse
therefore gives an algorithm for computing expansions in our basis.
Moreover, the expansions in \cref{thm:mainA,thm:mainB} are obtained from those in \cref{thm:mainF}.
Consequently, we also obtain an algorithm for expressing finite multiple zeta values and symmetric multiple zeta values as $\bbQ$-linear combinations of values with Zhao indices of the same weight.
Furthermore, combining the expansion algorithm for \cref{thm:mainC} with the explicit stuffle identity in \cite[Proposition~6.5]{Yasuda2016}, we also obtain an algorithm for expressing any multiple zeta value as a $\bbQ$-linear combination of $\zeta^{\sh}_{\calS}$-values of the same weight, as in Yasuda's theorem \cref{thm:Yasuda}.

For every nonempty index $\bk$, write the expansion provided by \cref{thm:mainF}, modulo $\zeta^{\rcf}(2)$, as
\begin{equation}\label{eq:unique-expansion}
\zeta^{\rcf}_{\calS}(\bk)\equiv\sum_{\bl\in\bbI^{\rZ}_{\lvert\bk\rvert}}c_{\bk,\bl}\zeta^{\rcf}_{\calS}(\bl)\pmod{\zeta^{\rcf}(2)\calZ^{\rcf}_{\lvert\bk\rvert-2}}.
\end{equation}
Let $S_k$ denote the basis given by \cref{thm:mainF} for $k\geq0$, and set $S_k=\varnothing$ for $k<0$.
For every $k\geq1$, we have the disjoint union decomposition
\[
S_k=\{\zeta^{\rcf}_{\calS}(\bk) : \bk\in\bbI^{\rZ}_k\}\cup\zeta^{\rcf}(2)S_{k-2},
\]
where $\zeta^{\rcf}(2)S_{k-2}$ denotes the set obtained by multiplying each element of $S_{k-2}$ by $\zeta^{\rcf}(2)$.
It follows that
\[
\calZ^{\rcf}_k=\Span_{\bbQ}\{\zeta^{\rcf}_{\calS}(\bk) : \bk\in\bbI^{\rZ}_k\}\oplus\zeta^{\rcf}(2)\calZ^{\rcf}_{k-2}.
\]
Hence the elements of the set $\{\zeta^{\rcf}_{\calS}(\bk)\bmod{\zeta^{\rcf}(2)\calZ^{\rcf}_{k-2}} : \bk\in\bbI^{\rZ}_k\}$ form a basis of $\calZ^{\rcf}_k/\zeta^{\rcf}(2)\calZ^{\rcf}_{k-2}$.
In particular, the coefficients $c_{\bk,\bl}\in\bbQ$ are uniquely determined.

In contrast to the Hoffman expansion in \cref{thm:Hirose}, these coefficients need not belong to $\bbZ_{(2)}$.
Indeed,
\[
\zeta^{\rcf}_{\calS}(2,1,3)\equiv\frac{1}{2}\zeta^{\rcf}_{\calS}(3,2,1)\pmod{\zeta^{\rcf}(2)\calZ^{\rcf}_4}.
\]
Exact computations in weights at most $13$ verify that $2^{\iota_3(\bl)}c_{\bk,\bl}\in\bbZ_{(2)}$, where $\iota_3(\bl)$ denotes the number of consecutive $3$'s at the beginning of $\bl$.
However, this property does not hold in general.
For example,
\[
c_{(3,3,3,1,2,2),\,(2,2,3,2,2,2,1)}=
\frac{882580110056929867097}{2\cdot 1648980856616730273}.
\]
Since $\iota_3(2,2,3,2,2,2,1)=0$ and the displayed coefficient has $2$-adic valuation $-1$, this gives a counterexample in weight $14$.

For each positive integer $k$, we use the coefficients in \eqref{eq:unique-expansion} to define
\begin{equation}\label{eq:def_R_k}
\calR_k\coloneqq\Span_{\bbQ}\Biggl\{\sfw(\bk)-\sum_{\bl\in\bbI^{\rZ}_{\lvert\bk\rvert}}c_{\bk,\bl}\sfw(\bl) : \bk\in\bbI_k\setminus\bbI^{\rZ}_k\Biggr\}.
\end{equation}
Write $r_{\bk}\coloneqq\sfw(\bk)-\sum_{\bl\in\bbI^{\rZ}_{\lvert\bk\rvert}}c_{\bk,\bl}\sfw(\bl)$.
If $\bk,\bk'\in\bbI_k\setminus\bbI^{\rZ}_k$ are distinct, then the coefficients of $\sfw(\bk)$ and $\sfw(\bk')$ in $r_{\bk}$ are $1$ and $0$, respectively, while those in $r_{\bk'}$ are $0$ and $1$, respectively.
This observation shows that the elements of the set $\{r_{\bk} : \bk\in\bbI_k\setminus\bbI^{\rZ}_k\}$ form a basis of $\calR_k$, and hence
\begin{equation}\label{eq:dim_of_R_k}
\dim_{\bbQ}\calR_k=\#(\bbI_k\setminus\bbI^{\rZ}_k)=2^{k-1}-d_{k-3}.
\end{equation}

In our proof of \cref{thm:mainF}, we use Hirose's upper bound $\dim_{\bbQ}\calZ^{\rcf}_k\leq d_k$, so it remains to prove that the proposed basis elements are linearly independent over $\bbQ$.
The strategy is similar to that of Brown \cite{Brown2012}: linear independence is reduced to the invertibility of a certain matrix, which is ultimately established by $2$-adic computations.
The use of the following formula of Zagier is another common feature.
\begin{theorem}[Zagier \cite{Zagier2012}]\label{thm:Zagier232}
For nonnegative integers $a$ and $b$, we have
\[
\zeta(\{2\}^a,3,\{2\}^b)=\sum_{l=1}^{a+b+1}c_{a,b}^l\zeta(\{2\}^{a+b+1-l})\zeta(2l+1),
\]
where
\[
c_{a,b}^l=2(-1)^l\left(\binom{2l}{2a+2}-(1-4^{-l})\binom{2l}{2b+1}\right).
\]
\end{theorem}
Brown \cite[Section~9]{Brown2012} observed that his argument could be dualized and formulated in terms of Ihara's action rather than Goncharov's coproduct formula.
In line with this viewpoint, we adopt a formulation based on the Ihara product, so our argument differs in presentation from Brown's original proof.
Despite these similarities, as emphasized above, we do not use any motivic constructions.
Whereas Brown needs to lift Zagier's formula to an identity among motivic multiple zeta values, we use the original identity among real multiple zeta values directly.
We also use the alternating block notation of Charlton \cite{Charlton2021} where appropriate.
Since we work with Zhao indices rather than Hoffman indices, it is essential to introduce a matrix adapted to this family.

A related non-motivic adaptation of Brown's method was given by Willwacher \cite{WillwacherFree2026} in his proof of a freeness result for the Grothendieck--Teichm\"uller Lie algebra.
His argument also uses the Ihara action and $2$-adic estimates for a pairing matrix.

The matrix is defined as follows.
For $N\geq r\geq 1$, let $M_{N,r}$ be the square matrix of size $\binom{N-1}{r-1}\times\binom{N-1}{r-1}$ given by
\[
M_{N,r}\coloneqq(\xi_{\bk,\bl})_{(\bk,\bl)\in\bbI_{N,r}\times\bbI_{N,r}},
\]
where, for $\bl=(l_1,\dots,l_r)$, the entries are defined by
\[
\xi_{\bk,\bl}\coloneqq\lambda_{\bl}\cdot[v_1^{2l_1+1}\cdots v_r^{2l_r+1}]\left([\sfw(\sfz(\bk))e_1](\Phi_{v_1,\dots,v_r}^{-1}e_1\Phi_{v_1,\dots,v_r})\right).
\]
Here, for each positive integer $l$, we set
\begin{equation}\label{eq:def_of_lambda_l}
\lambda_l\coloneqq\frac{4^{l-1}}{2l\cdot\zeta(2l+1)},
\end{equation}
and write $\lambda_{\bl}=\lambda_{l_1}\cdots\lambda_{l_r}$.
The variables $v_1,\dots,v_r$ are commuting indeterminates, and the series $\Phi_{v_1,\dots,v_r}\in\bbR[v_1,\dots,v_r]\fps{e_0,e_1}$ is defined in \cref{def_of_mainassociators}.
A priori, this matrix has real entries, but we will show that all its entries lie in $\bbZ_{(2)}$.
\begin{maintheorem}\label{thm:mainG}
Let $N\geq r\geq1$.
Then $M_{N,r}$ is invertible over $\bbQ$.
\end{maintheorem}
The matrix $M_{N,r}$ arises by evaluating the proposed formal basis elements at a family of associators and extracting coefficients in the parameters.
Its invertibility allows us to separate the proposed basis elements and prove their linear independence.
Combined with Hirose's upper bound, this yields \cref{thm:mainF}.
\section{The Ihara product and associators}
In this section, we fix our conventions for the Ihara product and introduce the associators used to define evaluation maps on $\calZ^{\rcf}$.
We then establish the coefficient formulas and vanishing properties needed for the matrix computations and the proof of linear independence.
\subsection{Ihara product}
Let $K$ be a field of characteristic $0$, and let $K\fps{e_0,e_1}$ denote the ring of noncommutative formal power series in the two variables $e_0$ and $e_1$ with coefficients in $K$.
For a word $w\in\calH$ and a series $F\in K\fps{e_0,e_1}$, we denote the coefficient of $w$ in $F$ by $[w]F$.
The usual multiplication is induced by concatenation of words.
A series $F$ with $[1]F=1$ is invertible for this multiplication, and $F^{-1}$ will always denote this inverse.
The empty word is denoted by $1$.
For a word $w$, let $\overleftarrow{w}$ denote its reversal, and write $[w]F$ for its coefficient in a series $F$.
In particular, $[1]F$ is the constant term of $F$.
We extend the coefficient notation linearly to elements of $\calH$.

A series $F\in K\fps{e_0,e_1}$ is called \emph{group-like} if $[1]F=1$ and $[w_1\sh w_2]F=[w_1]F[w_2]F$ for all words $w_1,w_2$.
If $F$ is group-like, then $[w]F^{-1}=(-1)^{\lvert w\rvert}[\overleftarrow{w}]F$ for every word $w$.
This follows from the antipode identity $\sum_{w=w_1w_2}(-1)^{\lvert w_1\rvert}\overleftarrow{w_1}\sh w_2=0$, which holds for every nonempty word $w$.
We say that $F\in K\fps{e_0,e_1}$ is \emph{commutator group-like} if it is group-like and satisfies $[e_0]F=[e_1]F=0$.

Define the maps $\id$ and $\sigma$ on $K\fps{e_0,e_1}$ by $\id(F)=F$ and $\sigma(F)=F(-e_0,-e_1)$.
For $F\in K\fps{e_0,e_1}$ with $[1]F=1$, define the operator $T_F$ by $T_F(G)=G(e_0,F^{-1}e_1F)$ for $G\in K\fps{e_0,e_1}$.
This substitution is well-defined because both substituted series have zero constant term.
The map $T_F$ is a $K$-algebra homomorphism with respect to the usual multiplication.
The \emph{Ihara product} $\circ$ is defined by
\[
F\circ G\coloneqq F\cdot T_F(G)=F\cdot G(e_0,F^{-1}e_1F),
\]
where $\cdot$ denotes the usual multiplication of noncommutative formal power series.
The Ihara product is associative on the set of series with constant term $1$, with identity element $1$.
For $F,G\in K\fps{e_0,e_1}$ with $[1]F=[1]G=1$, it is straightforward to check that $T_{F\circ G}=T_FT_G$, where $T_FT_G$ denotes composition, with $T_G$ applied first.
\subsection{Associators}
We construct a family of formal power series satisfying the pentagon equation and use Furusho's theorem to obtain evaluation maps on $\calZ^{\rcf}$.
We also express the evaluations of formal symmetric multiple zeta values as coefficients of the corresponding adjoint series.

Let $\Phi_{\rKZ}\in\bbR\fps{e_0,e_1}$ be the \emph{Drinfel'd associator}, also known as the \emph{KZ associator}, where KZ stands for Knizhnik--Zamolodchikov.
We use the normalization in which
\[
[\overleftarrow{\sfw(\bk)}]\Phi_{\rKZ}=(-1)^{\dep(\bk)}\zeta^{\sh}(\bk)
\]
for every index $\bk$.
This series is commutator group-like and satisfies the pentagon equation; see \cite{Drinfeld1991,Furusho2007}.

The formal power series $\Phi^-_{\rKZ}$ with constant term $1$, uniquely characterized by
\begin{equation}\label{eq:def_of_phi^-}
\Phi_{\rKZ}=\Phi^-_{\rKZ}\circ\sigma(\Phi_{\rKZ}),
\end{equation}
plays an important role in this paper.
We follow Furusho's notation; see \cite[Lemma~2.25]{Furusho2007}.
After translating the conventions for word order and signs, this is the series called \emph{Deligne's associator} by Brown.
Its coefficients encode the regularized values at $1$ of single-valued multiple polylogarithms; see \cite[Section~2.2]{Furusho2007} and \cite[Corollary~5.4]{Brown2014}.

Furusho proved that $\Phi^-_{\rKZ}\in\GRT_1(\bbR)$ in \cite[Lemma~3.3 and Remark~3.4]{Furusho2007}, where $\GRT_1(\bbR)$ denotes the group of real points of the prounipotent graded Grothendieck--Teichm\"uller group, with group law given by the Ihara product.
\begin{definition}\label{def_of_mainassociators}
Let $d$ be a nonnegative integer, and let $s,v_1,\dots,v_d$ be commuting variables.
Define
\[
\Phi_{v_1,\dots, v_d}\coloneqq \Phi^-_{\rKZ}(v_1e_0, v_1e_1)\circ\cdots\circ\Phi^-_{\rKZ}(v_de_0, v_de_1)
\]
and
\[
\Phi_{s;v_1,\dots,v_d}\coloneqq\Phi_{v_1,\dots,v_d}\circ\Phi_{\rKZ}(se_0,se_1).
\]
We set $\Phi_{\varnothing}=1$, interpreting an empty Ihara product as $1$. Every coefficient of these series is a polynomial in the indicated parameters with real coefficients. In particular,
\[
\Phi_{s;v_1,\dots,v_d}\in\bbR[s,v_1,\dots,v_d]\fps{e_0,e_1}.
\]
Then $\Phi_{s;\varnothing}=\Phi_{\rKZ}(se_0,se_1)$ and $\Phi_{0;v_1,\dots,v_d}=\Phi_{v_1,\dots,v_d}$.
\end{definition}
Using the identity $T_{F\circ G}=T_FT_G$, we obtain
\begin{equation}\label{eq:adj-Phi=T...T(e_1)}
\Phi_{v_1,\dots,v_d}^{-1}e_1\Phi_{v_1,\dots,v_d}=T_{\Phi_{v_1,\dots,v_d}}(e_1)=T_{\Phi^-_{\rKZ}(v_1e_0,v_1e_1)}\cdots T_{\Phi^-_{\rKZ}(v_de_0,v_de_1)}(e_1).
\end{equation}
\begin{lemma}\label{lem:pentagon_for_Phi_sv}
The formal power series $\Phi_{s;v_1,\dots,v_d}$ is commutator group-like and satisfies the pentagon equation.
\end{lemma}
\begin{proof}
For a field $K$ of characteristic $0$, let $\GRT_1(K)$ denote the group of degenerate associators with parameter $0$, with group law given by the Ihara product, and let $M_{\mu}(K)$ denote the set of associators with parameter $\mu\in K^{\times}$.
By \cite[Section~5, Proposition~5.5 and the preceding discussion]{Drinfeld1991}, the group $\GRT_1(K)$ acts on $M_{\mu}(K)$ from the left.
This means that $F\circ G\in M_{\mu}(K)$ whenever $F\in\GRT_1(K)$ and $G\in M_{\mu}(K)$.
Note that our multiplication is opposite to Drinfel'd's convention, so his right action becomes a left action in our convention.

It suffices to work over $K=\bbC(s,v_1,\dots,v_d)$.
Since $\Phi_{\rKZ}\in M_{2\pi\sqrt{-1}}(K)$ and $\Phi^-_{\rKZ}\in\GRT_1(K)$, simultaneous rescaling of the variables gives $\Phi_{\rKZ}(se_0,se_1)\in M_{2\pi\sqrt{-1}s}(K)$ and $\Phi^-_{\rKZ}(v_ie_0,v_ie_1)\in\GRT_1(K)$ for each $i$.
Applying the action successively, we conclude that $\Phi_{s;v_1,\dots,v_d}\in M_{2\pi\sqrt{-1}s}(K)$.
In particular, $\Phi_{s;v_1,\dots,v_d}$ is commutator group-like and satisfies the pentagon equation.
\end{proof}
Since every coefficient of $\Phi_{s;v_1,\dots,v_d}$ belongs to $\bbR[s,v_1,\dots,v_d]$, the group-like identities, the vanishing of the linear terms, and the pentagon equation hold over this polynomial ring.
They therefore remain valid under specialization of any of the parameters to zero.

In this paper, we use Furusho's result \cite{Furusho2022} only in the direction from the pentagon equation to the confluence relations.
\begin{theorem}[Furusho {\cite[Theorem~4]{Furusho2022}}]\label{thm:pentagon=confluence}
Suppose that $\Phi\in K\fps{e_0,e_1}$ is commutator group-like and satisfies the pentagon equation.
Then there is a well-defined $\bbQ$-algebra homomorphism
\[
\langle\cdot,\Phi\rangle\colon\calZ^{\rcf}\longrightarrow K
\]
given by $\langle w\bmod{\calI_{\rCF}},\Phi\rangle=[\overleftarrow{w}]\Phi$ for $w\in\calH^0$.
\end{theorem}
The reversal appears because Furusho uses the opposite word order to the original convention of Hirose and Sato, which we follow here; see \cite[Section~2]{Furusho2022}.
Although the result in \cite{Furusho2022} is stated only for $K=\bbQ$, Furusho's proof applies to any field $K$ of characteristic $0$.
For $w\in\calH^0$, we also write $\langle w,\Phi\rangle$ for $\langle w\bmod{\calI_{\rCF}},\Phi\rangle$.
\begin{lemma}
Suppose that $\Phi\in K\fps{e_0,e_1}$ is commutator group-like and satisfies the pentagon equation.
Then, for every index $\bk$, not necessarily admissible, we have
\begin{equation}\label{eq:cf-Phi_general_index}
\langle\zeta^{\rcf}(\bk),\Phi\rangle=(-1)^{\dep(\bk)}[\overleftarrow{\sfw(\bk)}]\Phi
\end{equation}
and
\begin{equation}\label{eq:sym_cf-Phi}
\langle\zeta^{\rcf}_{\calS}(\bk),\Phi\rangle=(-1)^{\lvert\bk\rvert+\dep(\bk)}[\sfw(\bk)e_1](\Phi^{-1}e_1\Phi).
\end{equation}
\end{lemma}
\begin{proof}
The identity \eqref{eq:cf-Phi_general_index} follows from the definition of $\zeta^{\rcf}(\bk)$ in \eqref{eq:def_of_zeta^cf} and the fact that $\Phi$ is group-like and $[e_0]\Phi=[e_1]\Phi=0$.

Write $\bk=(k_1,\dots,k_r)$, and let $0\leq i\leq r$.
If we match the central $e_1$ in $\Phi^{-1}e_1\Phi$ with the $(i+1)$st occurrence of $e_1$ in $\sfw(\bk)e_1$, then the subwords coming from $\Phi^{-1}$ and $\Phi$ are $\sfw(k_1,\dots,k_i)$ and $\overleftarrow{\sfw(k_r,\dots,k_{i+1})}$, respectively.
Since $[\sfw(k_1,\dots,k_i)]\Phi^{-1}=(-1)^{k_1+\cdots+k_i}[\overleftarrow{\sfw(k_1,\dots,k_i)}]\Phi$, the contribution of this decomposition to $[\sfw(\bk)e_1](\Phi^{-1}e_1\Phi)$ is, by \eqref{eq:cf-Phi_general_index},
\begin{align*}
&(-1)^{k_1+\cdots+k_i}[\overleftarrow{\sfw(k_1,\dots,k_i)}]\Phi\cdot[\overleftarrow{\sfw(k_r,\dots,k_{i+1})}]\Phi\\
&=(-1)^{k_1+\cdots+k_i+r}\langle\zeta^{\rcf}(k_1,\dots,k_i),\Phi\rangle\langle\zeta^{\rcf}(k_r,\dots,k_{i+1}),\Phi\rangle\\
&=\langle(-1)^{k_1+\cdots+k_i+r}\zeta^{\rcf}(k_1,\dots,k_i)\zeta^{\rcf}(k_r,\dots,k_{i+1}),\Phi\rangle.
\end{align*}
The desired identity now follows by comparison with the defining formula \eqref{eq:def_of_sym-zeta^cf}.
\end{proof}
\subsection{Coefficient formulas for the Ihara product}
Although well known, we first record coefficient formulas for the Ihara product and the substitution operator $T_F$.
When $F$ is group-like, the first formula expresses the duality between the Ihara product and the Goncharov coproduct on the shuffle algebra of words; see \cite[Lemma~3.25 and its proof]{BurmesterConfuriusKuhn}.
We then use these formulas to establish vanishing properties of $\Phi^-_{\rKZ}$ and to control the block degree in iterated applications of $T_{\Phi^-_{\rKZ}(ve_0,ve_1)}-\id$.
\begin{lemma}\label{lem:circ_calc}
Set $\epsilon(e_0)\coloneqq0$ and $\epsilon(e_1)\coloneqq1$.
Let $F,G\in K\fps{e_0,e_1}$ with $[1]F=1$.
Then, for every word $w$, we have
\[
[w](F\circ G)=\sum_{k=0}^{\lvert w\rvert}\sum_{\substack{w=w_0u_1w_1\cdots u_kw_k\\u_1,\dots,u_k\in\{e_0,e_1\}}}[u_1\cdots u_k]G \, \prod_{i=0}^k[w_i]F^{\epsilon(u_i)-\epsilon(u_{i+1})},
\]
where $u_0=e_1$ and $u_{k+1}=e_0$, and each word $w_i$ in the decomposition $w=w_0u_1w_1\cdots u_kw_k$ is allowed to be empty.
With the convention $u_0=u_{k+1}=e_0$ instead, we also have
\[
[w]T_F(G)=\sum_{k=0}^{\lvert w\rvert}\sum_{\substack{w=w_0u_1w_1\cdots u_kw_k\\u_1,\dots,u_k\in\{e_0,e_1\}}}[u_1\cdots u_k]G \, \prod_{i=0}^k[w_i]F^{\epsilon(u_i)-\epsilon(u_{i+1})}.
\]
\end{lemma}
\begin{proof}
First, let $k\geq1$ and fix $u_1,\dots,u_k\in\{e_0,e_1\}$.
With $u_0=e_1$ and $u_{k+1}=e_0$, we have
\begin{align*}
F\circ(u_1\cdots u_k)&=F(F^{-\epsilon(u_1)}u_1F^{\epsilon(u_1)})\cdots(F^{-\epsilon(u_k)}u_kF^{\epsilon(u_k)})\\
&=F^{\epsilon(u_0)-\epsilon(u_1)}u_1F^{\epsilon(u_1)-\epsilon(u_2)}u_2\cdots u_kF^{\epsilon(u_k)-\epsilon(u_{k+1})}.
\end{align*}
Expanding each factor as
\[
F^{\epsilon(u_i)-\epsilon(u_{i+1})}=\sum_{w_i}([w_i]F^{\epsilon(u_i)-\epsilon(u_{i+1})})w_i,
\]
where $w_i$ runs over all words, including the empty word, gives
\[
[w](F\circ(u_1\cdots u_k))=\sum_{w=w_0u_1w_1\cdots u_kw_k}\prod_{i=0}^k[w_i]F^{\epsilon(u_i)-\epsilon(u_{i+1})}.
\]
The case $k=0$ is given by $F\circ1=F$.
Since
\[
G=\sum_{k=0}^{\infty}\sum_{u_1,\dots,u_k\in\{e_0,e_1\}}([u_1\cdots u_k]G)u_1\cdots u_k
\]
and the map $G\mapsto F\circ G$ is $K$-linear, summing the preceding coefficient formula gives the first identity.
Only $k\leq\lvert w\rvert$ can contribute, since every word occurring in $F\circ(u_1\cdots u_k)$ has length at least $k$.

For the second identity, omit the initial factor $F$ in the first calculation.
The same argument then applies with $u_0=u_{k+1}=e_0$.
The case $k=0$ is given by $T_F(1)=1$.
\end{proof}
Although \cref{lem:(xy)^n=0} and \cref{cor:degree2vanish} are known, we include proofs here, as both results follow readily from \cref{lem:circ_calc}.
A nonempty word is called \emph{alternating} if any two consecutive letters are distinct.
\begin{lemma}\label{lem:(xy)^n=0}
For every positive integer $n$, we have
\[
[(e_0e_1)^n]\Phi^-_{\rKZ}=0.
\]
\end{lemma}
\begin{proof}
Since $\Phi^-_{\rKZ}$ is group-like, $[(e_1e_0)^j](\Phi^-_{\rKZ})^{-1}=[(e_0e_1)^j]\Phi^-_{\rKZ}$.
Applying \cref{lem:circ_calc} to \eqref{eq:def_of_phi^-}, we obtain
\[
[(e_0e_1)^n]\Phi_{\rKZ}=\sum_{k=0}^{2n}\sum_{\substack{(e_0e_1)^n=w_0u_1w_1\cdots u_kw_k \\ u_1,\dots, u_k\in\{e_0, e_1\}}}[u_1\cdots u_k]\Phi_{\rKZ}(-e_0, -e_1)\prod_{i=0}^k[w_i](\Phi^-_{\rKZ})^{\epsilon(u_i)-\epsilon(u_{i+1})}.
\]
Here $u_0=e_1$ and $u_{k+1}=e_0$.
We proceed by induction on $n$.
Assume that $[(e_0e_1)^j]\Phi^{-}_{\rKZ}=0$ for all $1\leq j<n$.
Let $1\leq k\leq 2n-1$, and fix a decomposition $(e_0e_1)^n=w_0u_1w_1\cdots u_kw_k$.
Since $\sum_{i=0}^k\lvert w_i\rvert=2n-k>0$, at least one of the words $w_i$ is nonempty.
Fix such an $i$.
We now show that $[w_i](\Phi^-_{\rKZ})^{\epsilon(u_i)-\epsilon(u_{i+1})}=0$.
If $u_i=u_{i+1}$, then
\[
[w_i](\Phi^-_{\rKZ})^{\epsilon(u_i)-\epsilon(u_{i+1})}=[w_i](\Phi^-_{\rKZ})^0=[w_i]1=0.
\]
Suppose that $u_i=e_1$ and $u_{i+1}=e_0$.
Since $e_1w_ie_0$ is an alternating word, we have $w_i=(e_0e_1)^j$ for some $j\geq 1$.
The bound $2j=\lvert w_i\rvert\leq 2n-k\leq 2n-1$ implies that $1\leq j<n$.
Thus, by the induction hypothesis,
\[
[w_i](\Phi^-_{\rKZ})^{\epsilon(u_i)-\epsilon(u_{i+1})}=[(e_0e_1)^j]\Phi^-_{\rKZ}=0.
\]
If $u_i=e_0$ and $u_{i+1}=e_1$, then $e_0w_ie_1$ is an alternating word, so $w_i=(e_1e_0)^j$ for some $1\leq j<n$.
Hence
\[
[w_i]\bigl(\Phi^-_{\rKZ}\bigr)^{\epsilon(u_i)-\epsilon(u_{i+1})}=
[(e_1e_0)^j]\bigl(\Phi^-_{\rKZ}\bigr)^{-1}=[(e_0e_1)^j]\Phi^-_{\rKZ}=0.
\]
For $k=0$, we have $w_0=(e_0e_1)^n$, and the corresponding term is $[(e_0e_1)^n]\Phi^-_{\rKZ}$.
For $k=2n$, we have $u_1\cdots u_{2n}=(e_0e_1)^n$, and the corresponding term is $[(e_0e_1)^n]\Phi_{\rKZ}(-e_0,-e_1)$.
Combining these observations, we obtain
\[
[(e_0e_1)^n]\Phi_{\rKZ}=[(e_0e_1)^n]\Phi^-_{\rKZ}+[(e_0e_1)^n]\Phi_{\rKZ}(-e_0,-e_1).
\]
Since
\[
[(e_0e_1)^n]\Phi_{\rKZ}=[(e_0e_1)^n]\Phi_{\rKZ}(-e_0,-e_1),
\]
we conclude that $[(e_0e_1)^n]\Phi^-_{\rKZ}=0$.

When $n=1$, the only intermediate terms are those with $k=1$.
There are two possible decompositions $e_0e_1=w_0u_1w_1$.
If $w_0=1$, $u_1=e_0$, and $w_1=e_1$, the corresponding product contains the vanishing factor $[e_1]\bigl(\Phi^-_{\rKZ}\bigr)^{0-0}=0$.
If $w_0=e_0$, $u_1=e_1$, and $w_1=1$, it contains the vanishing factor $[e_0]\bigl(\Phi^-_{\rKZ}\bigr)^{1-1}=0$.
Thus, the same argument as in the induction step establishes the base case $[e_0e_1]\Phi^-_{\rKZ}=0$.
\end{proof}
\begin{corollary}\label{cor:degree2vanish}
All coefficients of words of length $2$ in $\Phi^-_{\rKZ}$ vanish.
\end{corollary}
\begin{proof}
The case $n=1$ of the preceding lemma gives $[e_0e_1]\Phi^-_{\rKZ}=0$.
Since $\Phi^-_{\rKZ}$ is commutator group-like, we have $[e_0]\Phi^-_{\rKZ}=[e_1]\Phi^-_{\rKZ}=0$.
The shuffle identities $e_i\sh e_i=2e_i^2$ for $i=0,1$ and $e_0\sh e_1=e_0e_1+e_1e_0$ therefore imply that $[e_i^2]\Phi^-_{\rKZ}=0$ for $i=0,1$ and $[e_0e_1]\Phi^-_{\rKZ}+[e_1e_0]\Phi^-_{\rKZ}=0$.
Hence $[e_1e_0]\Phi^-_{\rKZ}=0$, completing the proof.
\end{proof}
For a word $w=u_1\cdots u_n$ with $u_i\in\{e_0,e_1\}$, we define its \emph{block degree} by
\[
\deg_{\rbl}(w)\coloneqq\#\{i : 1\leq i<n,\ u_i=u_{i+1}\}.
\]
See also \cite[Definition~3.1]{Keilthy2022}.
Thus, a nonempty word is alternating if and only if its block degree is $0$. We set $\deg_{\rbl}(1)=0$.
\begin{lemma}\label{lem:T-id}
Let $r$ be a positive integer, let $v_1,\dots,v_r$ be commuting indeterminates, and let $w$ be a word in $e_0,e_1$.
If $\deg_{\rbl}(e_0we_0)<r$, then
\begin{equation}\label{eq:T-id_zero}
[w](T_{\Phi^-_{\rKZ}(v_1e_0, v_1e_1)}-\id)\cdots(T_{\Phi^-_{\rKZ}(v_re_0, v_re_1)}-\id)(e_1)=0.
\end{equation}
If $\deg_{\rbl}(e_0we_0)=r$, then
\begin{equation}\label{eq:T-id_rec}
\begin{split}
&[w](T_{\Phi^-_{\rKZ}(v_1e_0, v_1e_1)}-\id)\cdots(T_{\Phi^-_{\rKZ}(v_re_0, v_re_1)}-\id)(e_1)\\
&=\sum_{k=1}^{\lvert w\rvert}\sum_{\substack{w=w_0u_1w_1\cdots u_kw_k \\ u_1,\dots, u_k\in\{e_0, e_1\} \\ \#\{i : 0\leq i\leq k, w_i\neq 1\}=1 \\ \deg_{\rbl}(u_0u_1\cdots u_{k+1})=r-1}}\left(\prod_{i=0}^k[w_i](\Phi^-_{\rKZ}(v_1e_0, v_1e_1))^{\epsilon(u_i)-\epsilon(u_{i+1})}\right)\\
&\quad\cdot[u_1\cdots u_k](T_{\Phi^-_{\rKZ}(v_2e_0, v_2e_1)}-\id)\cdots(T_{\Phi^-_{\rKZ}(v_re_0, v_re_1)}-\id)(e_1).
\end{split}
\end{equation}
In \eqref{eq:T-id_rec}, we set $u_0=u_{k+1}=e_0$, and each word $w_i$ is allowed to be empty. When $r=1$, the empty product of operators on the right-hand side is interpreted as $\id$.
\end{lemma}
\begin{proof}
Throughout the proof, we use the convention $u_0=u_{k+1}=e_0$.
In the second formula of \cref{lem:circ_calc}, the term for which all $w_i$ are empty is $[w]G$.
Subtracting this term gives
\begin{equation}\label{eq:T-id_fund}
[w](T_{\Phi^-_{\rKZ}}-\id)(G)=\sum_{k=1}^{\lvert w\rvert}\sum_{\substack{w=w_0u_1w_1\cdots u_kw_k \\ u_1,\dots, u_k\in\{e_0, e_1\} \\ \#\{i : 0\leq i\leq k, w_i\neq 1\}\geq 1}}[u_1\cdots u_k] \, G \prod_{i=0}^k[w_i](\Phi^-_{\rKZ})^{\epsilon(u_i)-\epsilon(u_{i+1})}.
\end{equation}
In the above expansion, consider a decomposition $w=w_0u_1w_1\cdots u_kw_k$ such that
\[
\prod_{i=0}^k[w_i](\Phi^-_{\rKZ})^{\epsilon(u_i)-\epsilon(u_{i+1})}\neq 0.
\]
We claim that
\begin{equation}\label{eq:diff_blockdeg_inequality}
\deg_{\rbl}(e_0we_0)-\deg_{\rbl}(u_0u_1\cdots u_{k+1})\geq \#\{i : w_i\neq 1\}.
\end{equation}
First, we show that $\deg_{\rbl}(u_iw_iu_{i+1})\geq 1$ whenever $w_i\neq 1$.
Fix such an $i$.
If $u_i=u_{i+1}$, then
\[
[w_i](\Phi^-_{\rKZ})^{\epsilon(u_i)-\epsilon(u_{i+1})}=[w_i]1=0,
\]
contrary to our assumption.
Thus $u_i\neq u_{i+1}$.
Suppose, for a contradiction, that $u_iw_iu_{i+1}$ is an alternating word.
If $u_i=e_1$ and $u_{i+1}=e_0$, then, since $w_i\neq 1$, we have $w_i=(e_0e_1)^j$ for some $j\geq 1$.
By \cref{lem:(xy)^n=0}, we obtain
\[
[w_i](\Phi^-_{\rKZ})^{1-0}=[(e_0e_1)^j]\Phi^-_{\rKZ}=0,
\]
a contradiction.
If $u_i=e_0$ and $u_{i+1}=e_1$, then similarly $w_i=(e_1e_0)^j$ for some $j\geq 1$.
By \cref{lem:(xy)^n=0}, we obtain
\[
[w_i](\Phi^-_{\rKZ})^{0-1}=[(e_1e_0)^j](\Phi^-_{\rKZ})^{-1}=[(e_0e_1)^j]\Phi^-_{\rKZ}=0,
\]
again a contradiction.
This proves that $\deg_{\rbl}(u_iw_iu_{i+1})\geq 1$ whenever $w_i\neq 1$.
Next, we compute the block degrees.
We have
\[
\deg_{\rbl}(e_0we_0)=\sum_{i=0}^k \deg_{\rbl}(u_iw_iu_{i+1})=\sum_{\substack{0\leq i\leq k\\ w_i=1}}\deg_{\rbl}(u_iw_iu_{i+1})+\sum_{\substack{0\leq i\leq k\\ w_i\neq 1}}\deg_{\rbl}(u_iw_iu_{i+1})
\]
and
\begin{align*}
\deg_{\rbl}(u_0u_1\cdots u_{k+1})
&=\sum_{i=0}^k \deg_{\rbl}(u_iu_{i+1})=\sum_{\substack{0\leq i\leq k\\ w_i=1}}\deg_{\rbl}(u_iu_{i+1})+\sum_{\substack{0\leq i\leq k\\ w_i\neq 1}}\deg_{\rbl}(u_iu_{i+1})\\
&=\sum_{\substack{0\leq i\leq k\\ w_i=1}}\deg_{\rbl}(u_iw_iu_{i+1}).
\end{align*}
Here the last equality follows from the fact that $u_i\neq u_{i+1}$, and hence $\deg_{\rbl}(u_iu_{i+1})=0$, whenever $w_i\neq 1$.
Subtracting these identities, we obtain
\[
\deg_{\rbl}(e_0we_0)-\deg_{\rbl}(u_0u_1\cdots u_{k+1})=\sum_{\substack{0\leq i\leq k\\ w_i\neq 1}}\deg_{\rbl}(u_iw_iu_{i+1}).
\]
Since each summand on the right-hand side is at least $1$, the inequality \eqref{eq:diff_blockdeg_inequality} follows.

Since we are considering only the case $\#\{i : w_i\neq 1\}\geq 1$,
we obtain, in particular,
\[
\deg_{\rbl}(u_0u_1\cdots u_{k+1})
\leq \deg_{\rbl}(e_0we_0)-1.
\]
Let $v$ be a commuting indeterminate.
Since
\[
[w]\Phi^-_{\rKZ}(ve_0,ve_1)=v^{\lvert w\rvert}[w]\Phi^-_{\rKZ}, \quad
(\Phi^-_{\rKZ})^{-1}(ve_0,ve_1)=(\Phi^-_{\rKZ}(ve_0,ve_1))^{-1},
\]
the same inequality holds with $\Phi^-_{\rKZ}$ replaced by $\Phi^-_{\rKZ}(ve_0,ve_1)$.

For each $1\leq i\leq r$, set
\[
G_i\coloneqq(T_{\Phi^-_{\rKZ}(v_ie_0,v_ie_1)}-\id)\cdots(T_{\Phi^-_{\rKZ}(v_re_0,v_re_1)}-\id)(e_1).
\]
We also set $G_{r+1}\coloneqq e_1$.
We show that $[w]G_1\neq 0$ implies $\deg_{\rbl}(e_0we_0)\geq r$.
Suppose that $[w]G_1\neq 0$.
Since
\[
[w]G_1=[w](T_{\Phi^-_{\rKZ}(v_1e_0,v_1e_1)}-\id)(G_2)\neq 0,
\]
there exists a scattered subword $w^{(1)}$ of $w$ such that
\[
[w^{(1)}]G_2\neq 0,\qquad\deg_{\rbl}(e_0w^{(1)}e_0)\leq \deg_{\rbl}(e_0we_0)-1.
\]
If $r\geq 2$, we can apply the same argument to
\[
[w^{(1)}]G_2=[w^{(1)}]\bigl(T_{\Phi^-_{\rKZ}(v_2e_0,v_2e_1)}-\id\bigr)(G_3)\neq 0
\]
to obtain a scattered subword $w^{(2)}$ of $w^{(1)}$ such that
\[
[w^{(2)}]G_3\neq 0,\qquad\deg_{\rbl}(e_0w^{(2)}e_0)\leq \deg_{\rbl}(e_0w^{(1)}e_0)-1.
\]
Continuing in this way, we eventually obtain a word $w^{(r)}$ such that $[w^{(r)}]e_1\neq 0$.
Then, we must have $w^{(r)}=e_1$.
Thus $\deg_{\rbl}(e_0w^{(r)}e_0)=0$, and hence
\[
0\leq \deg_{\rbl}(e_0w^{(r-1)}e_0)-1\leq \cdots\leq \deg_{\rbl}(e_0we_0)-r.
\]
This proves that $\deg_{\rbl}(e_0we_0)\geq r$ and hence \eqref{eq:T-id_zero}.

Now suppose that $\deg_{\rbl}(e_0we_0)=r$.
By the argument in the first part of the proof, each nonzero term in the expansion of
\[
[w]G_1=[w](T_{\Phi^-_{\rKZ}(v_1e_0,v_1e_1)}-\id)(G_2)
\]
given by \eqref{eq:T-id_fund} corresponds to a decrease in block degree of at least the number of nonempty words among $w_0,\dots,w_k$, and hence of at least $1$.
Thus
\[
[w]G_1=\sum_{k=1}^{\lvert w\rvert}\sum_{\substack{w=w_0u_1w_1\cdots u_kw_k\\
u_1,\dots,u_k\in\{e_0,e_1\}\\\#\{i : 0\leq i\leq k,\; w_i\neq 1\}\geq 1\\
\deg_{\rbl}(u_0u_1\cdots u_{k+1})\leq r-1}}\left(\prod_{i=0}^k[w_i](\Phi^-_{\rKZ}(v_1e_0,v_1e_1))^{\epsilon(u_i)-\epsilon(u_{i+1})}\right)\cdot [u_1\cdots u_k]G_2.
\]
Among these terms, those with $\deg_{\rbl}(u_0u_1\cdots u_{k+1})<r-1$
vanish by \eqref{eq:T-id_zero}.
For any nonzero term with $\deg_{\rbl}(u_0u_1\cdots u_{k+1})=r-1$, the block degree drops by exactly $1$, so we must have $\#\{i : 0\leq i\leq k,\; w_i\neq 1\}=1$.
This proves \eqref{eq:T-id_rec}.
\end{proof}
\section{An explicit formula for $M_{N,r}$ and its invertibility}
\subsection{An explicit formula for $M_{N,r}$}
We first compute the coefficients $\varphi_l(j)$ explicitly using Zagier's formula.
We then define the quantities $\eta_{\bl}$ recursively and show that their values give the entries of $M_{N,r}$.

Every nonempty word $w\in\calH$ admits a unique decomposition $w=B_1B_2\cdots B_r$ into nonempty alternating blocks, where each $B_i$ is an alternating word and the last letter of $B_i$ coincides with the first letter of $B_{i+1}$ for $1\leq i<r$; see \cite[Lemma~3.2]{Charlton2021}.
If the first letter of $B_1$ is $u$ and $\lvert B_i\rvert=b_i$, we write $\block(w)\coloneqq(u;b_1,\dots,b_r)$; see \cite[Definition~3.3]{Charlton2021}.
In this paper, however, the initial letter will always be clear from the context, so we omit it and write $\block(w)\coloneqq(b_1,\dots,b_r)$.

Let $l$ be a positive integer.
For $0\leq i\leq 2l+1$, let $w_{l,i}$ be the word determined by
\[
\block(e_1w_{l,i}e_0)=(i+1,2l+2-i).
\]
In particular, $\lvert w_{l,i}\rvert=2l+1$ and $\deg_{\rbl}(e_1w_{l,i}e_0)=1$.
Using the words $w_{l,i}$ and the series $\Phi^-_{\rKZ}$ defined in \eqref{eq:def_of_phi^-}, for each integer $j$, define $\varphi_l(j)\in\bbR$ by
\[
\varphi_l(j)\coloneqq\begin{cases}
\displaystyle\lambda_l\sum_{i=0}^j[w_{l,i}]\Phi^-_{\rKZ} & \text{if } 0\leq j\leq 2l,\\
0 & \text{if } j<0 \text{ or } j>2l,
\end{cases}
\]
where $\lambda_l$ is as defined in \eqref{eq:def_of_lambda_l}.
\begin{lemma}\label{lem:phi-fund}
For each $l, j$, we have
\begin{equation}\label{eq:phi-sym}
\varphi_l(j)=\varphi_l(2l-j),
\end{equation}
and
\begin{equation}\label{eq:phi-sum}
\sum_{j=0}^{2l}(-1)^j\varphi_l(j)=\frac{\varphi_l(0)}{2}.
\end{equation}
\end{lemma}
\begin{proof}
Assume $0\leq j\leq 2l$.
Let $\theta$ denote the automorphism of $\calH$ that interchanges $e_0$ and $e_1$.
The duality relation $\Phi^-_{\rKZ}(e_1,e_0)=(\Phi^-_{\rKZ}(e_0,e_1))^{-1}$ and the inverse-coefficient formula for group-like series imply that, for every word $w$ of odd length,
\[
[\theta(\overleftarrow{w})]\Phi^-_{\rKZ}=[\overleftarrow{w}](\Phi^-_{\rKZ})^{-1}=-[w]\Phi^-_{\rKZ}.
\]
Since $\theta(\overleftarrow{w_{l,i}})=w_{l,2l+1-i}$, we have
\[
[w_{l,2l+1-i}]\Phi^-_{\rKZ}=-[w_{l,i}]\Phi^-_{\rKZ}.
\]
This proves \eqref{eq:phi-sym}.
For $j<0$ or $j>2l$, both sides vanish by definition.

By
\[
e_0\sh(e_0e_1)^l=2\sum_{\substack{a, b\geq 0 \\ a+b=l-1}}(e_0e_1)^ae_0^2e_1(e_0e_1)^b+(e_0e_1)^le_0,
\]
$[e_0\sh(e_0e_1)^l]\Phi^-_{\rKZ}=[e_0]\Phi^-_{\rKZ}[(e_0e_1)^l]\Phi^-_{\rKZ}=0$ from \cref{lem:(xy)^n=0}, $(e_0e_1)^ae_0^2e_1(e_0e_1)^b=w_{l,2a+1}$, and $(e_0e_1)^le_0=w_{l,2l+1}$, we have
\[
2\sum_{a=0}^{l-1}(\varphi_l(2a+1)-\varphi_l(2a))-\varphi_l(2l)=0.
\]
Consequently,
\[
\sum_{j=0}^{2l}(-1)^j\varphi_l(j)=\varphi_l(2l)-\sum_{a=0}^{l-1}(\varphi_l(2a+1)-\varphi_l(2a))=\frac{\varphi_l(2l)}2=\frac{\varphi_l(0)}2.
\]
This proves \eqref{eq:phi-sum}.
\end{proof}
\begin{lemma}\label{lem:another_def_of_phi}
For every integer $j$, we have
\[
\varphi_l(j)=\lambda_l\sum_{\substack{0\leq i\leq 2l+1 \\ i\leq j}}[w_{l,i}]\Phi^-_{\rKZ}.
\]
\end{lemma}
\begin{proof}
By definition,
\[
\lambda_l\sum_{i=0}^{2l+1}[w_{l,i}]\Phi^-_{\rKZ}=\varphi_l(2l)+\lambda_l[w_{l,2l+1}]\Phi^-_{\rKZ}.
\]

By the proof of \cref{lem:phi-fund}, $\varphi_l(2l)=\varphi_l(0)$ and $[w_{l,2l+1}]\Phi^-_{\rKZ}=-[w_{l,0}]\Phi^-_{\rKZ}$.
Thus we have
\[
\lambda_l\sum_{i=0}^{2l+1}[w_{l,i}]\Phi^-_{\rKZ}=0.
\]
For $0\leq j\leq2l$, the formula is the definition of $\varphi_l(j)$.
For $j<0$, the sum is empty, and for $j>2l$, it is the sum of all the coefficients above, which is zero.
\end{proof}
\begin{lemma}\label{lem:phi-phi}
Let $l\geq1$ and $j,j'\geq0$ be integers satisfying $j+j'\geq2l$.
Let $w_1$ and $w_2$ be alternating words such that $\deg_{\rbl}(w_1w_2)=1$, $\lvert w_1\rvert=j+1$, and $\lvert w_2\rvert=j'+1$.
Then we have
\[
\lambda_l\sum_{\substack{w_1w_2=puwu'q \\ u, u'\in\{e_0, e_1\}, \ u\neq u' \\ \lvert w\rvert=2l+1 \\ \deg_{\rbl}(uwu')=1}}[w](\Phi^-_{\rKZ})^{\epsilon(u)-\epsilon(u')}=\varphi_l(j)-\varphi_l(j').
\]
Here $p$ and $q$ are possibly empty words.
\end{lemma}
\begin{proof}
For a decomposition $w_1w_2=puwu'q$, write $w=w_{\rL}w_{\rR}$ with $\lvert w_{\rL}\rvert=i$ and $\lvert w_{\rR}\rvert=2l+1-i$, so that $w_1=puw_{\rL}$ and $w_2=w_{\rR}u'q$.
Then $\lvert uw_{\rL}\rvert\leq \lvert w_1\rvert$ implies $i\leq j$, while $\lvert w_{\rR}u'\rvert\leq \lvert w_2\rvert$ implies $2l+1-i\leq j'$.
Thus $i$ must satisfy the three conditions $0\leq i\leq 2l+1$, $i\leq j$, and $2l+1-i\leq j'$.
Conversely, each $i$ satisfying these three conditions corresponds to exactly one decomposition appearing in the sum on the left-hand side.
For the decomposition corresponding to $i$, we have
$\block(uwu')=(i+1,2l+2-i)$.
If $u=e_1$ and $u'=e_0$, then $w=w_{l,i}$, and hence $[w](\Phi^-_{\rKZ})^{\epsilon(u)-\epsilon(u')}=[w_{l,i}]\Phi^-_{\rKZ}$.
If $u=e_0$ and $u'=e_1$, then $\theta(w)=w_{l,i}$, and the duality relation for $\Phi^-_{\rKZ}$ gives $[w](\Phi^-_{\rKZ})^{\epsilon(u)-\epsilon(u')}=[\theta(w)]\Phi^-_{\rKZ}=[w_{l,i}]\Phi^-_{\rKZ}$.
We therefore obtain
\[
\lambda_l\sum_{\substack{w_1w_2=puwu'q\\u,u'\in\{e_0,e_1\},\ u\neq u'\\
\lvert w\rvert=2l+1\\\deg_{\rbl}(uwu')=1}}[w](\Phi^-_{\rKZ})^{\epsilon(u)-\epsilon(u')}=\lambda_l\sum_{\substack{0\leq i\leq 2l+1\\i\leq j\\2l+1-i\leq j'}}[w_{l,i}]\Phi^-_{\rKZ}.
\]
Noting that $2l-j'\leq j$ and using \cref{lem:another_def_of_phi} together with \eqref{eq:phi-sym}, we can further compute
\begin{align*}
\lambda_l\sum_{\substack{0\leq i\leq 2l+1\\i\leq j\\2l+1-i\leq j'}}[w_{l,i}]\Phi^-_{\rKZ}&=\lambda_l\sum_{\substack{0\leq i\leq 2l+1\\i\leq j}}[w_{l,i}]\Phi^-_{\rKZ}-\lambda_l\sum_{\substack{0\leq i\leq 2l+1\\i\leq 2l-j'}}[w_{l,i}]\Phi^-_{\rKZ}\\&=\varphi_l(j)-\varphi_l(2l-j')=\varphi_l(j)-\varphi_l(j').\qedhere
\end{align*}
\end{proof}
\begin{proposition}
For every positive integer $l$ and every integer $j$ with $0\leq j\leq2l$, we have
\begin{equation}\label{eq:explicit_phi}
\varphi_l(j)=\frac{1-2\cdot 4^l-(-1)^j}{4l}\binom{2l}{j}.
\end{equation}
\end{proposition}
\begin{proof}
Let $a$ and $b$ be nonnegative integers, and set $n\coloneqq a+b+1$ and $w\coloneqq(e_0e_1)^be_0^2e_1(e_0e_1)^a$.
Then $\lvert w\rvert=2a+2b+3=2n+1$.
Since $\sfw(\{2\}^a,3,\{2\}^b)=\overleftarrow{w}$, we have $[w]\Phi_{\rKZ}=(-1)^n\zeta(\{2\}^a,3,\{2\}^b)$.
Moreover, $[w]\sigma(\Phi_{\rKZ})=-[w]\Phi_{\rKZ}$, and hence $[w](\Phi_{\rKZ}-\sigma(\Phi_{\rKZ}))=2(-1)^n\zeta(\{2\}^a,3,\{2\}^b)$.
Applying \cref{lem:circ_calc} to $\Phi_{\rKZ}=\Phi^-_{\rKZ}\circ\sigma(\Phi_{\rKZ})$, we obtain
\begin{equation}\label{eq:Phi-sigmaPhi}
[w](\Phi_{\rKZ}-\sigma(\Phi_{\rKZ}))=\sum_{k=0}^{\lvert w\rvert}\sum_{\substack{w=w_0u_1w_1\cdots u_kw_k \\ u_1,\dots,u_k\in\{e_0,e_1\} \\ \#\{i:w_i\neq 1\}\geq 1}}[u_1\cdots u_k]\sigma(\Phi_{\rKZ})\prod_{i=0}^k[w_i](\Phi^-_{\rKZ})^{\epsilon(u_i)-\epsilon(u_{i+1})},
\end{equation}
where $u_0=e_1$ and $u_{k+1}=e_0$.
Since $e_1we_0=(e_1e_0)^{b+1}\mid(e_0e_1)^{a+1}e_0$, we have $\block(e_1we_0)=(2b+2,2a+3)$ and $\deg_{\rbl}(e_1we_0)=1$.
The local block-degree argument in the proof of \cref{lem:T-id} applies equally with the boundary letters $u_0=e_1$ and $u_{k+1}=e_0$, including the case $k=0$.
For each nonzero term, it gives
\[
\deg_{\rbl}(e_1we_0)-\deg_{\rbl}(e_1u_1\cdots u_ke_0)\geq\#\{i:0\leq i\leq k,\ w_i\neq1\}.
\]
Since $\deg_{\rbl}(e_1we_0)=1$ and at least one $w_i$ is nonempty, $\deg_{\rbl}(e_1u_1\cdots u_ke_0)=0$ and exactly one $w_i$ is nonempty.
Write $\lvert w_i\rvert=2l+1$ for this unique nonempty $w_i$.
Since $e_1u_1\cdots u_ke_0$ is an alternating word, we have $u_1\cdots u_k=(e_0e_1)^{n-l}$.
We now express the right-hand side of \eqref{eq:Phi-sigmaPhi} as a sum over $l$.
When $l=0$, we have $w_i\in\{e_0,e_1\}$, so $[w_i]\Phi^-_{\rKZ}=[w_i](\Phi^-_{\rKZ})^{-1}=0$.
Thus, this case contributes nothing to the sum, and it suffices to consider $1\leq l\leq n$.
Since $[(e_0e_1)^{n-l}]\Phi_{\rKZ}=(-1)^{n-l}\zeta(\{2\}^{n-l})$, we obtain from \cref{lem:phi-phi} that
\begin{equation}\label{eq:abstract_Zagier}
2\zeta(\{2\}^a,3,\{2\}^b)=\sum_{l=1}^n\lambda_l^{-1}(-1)^l\zeta(\{2\}^{n-l})(\varphi_l(2b+1)-\varphi_l(2a+2)).
\end{equation}
We prove \eqref{eq:explicit_phi} by induction on $l$.
The case $l=1$ is covered by the same argument below, since there are simply no terms with smaller indices to eliminate.
Assume that \eqref{eq:explicit_phi} holds for all $1\leq h<l$.
Then
\[
\varphi_h(2b+1)=-\frac{4^h-1}{2h}\binom{2h}{2b+1},\quad \varphi_h(2a+2)=-\frac{4^h}{2h}\binom{2h}{2a+2}.
\]
Hence,
\[
\lambda_h^{-1}(-1)^h(\varphi_h(2b+1)-\varphi_h(2a+2))=2c^h_{a,b}\zeta(2h+1).
\]
Comparing Zagier's theorem \cref{thm:Zagier232} with \eqref{eq:abstract_Zagier} for $a+b+1=n=l$, we obtain
\[
\lambda_l^{-1}(-1)^l(\varphi_l(2b+1)-\varphi_l(2a+2))=2c^l_{a,b}\zeta(2l+1).
\]
By \eqref{eq:phi-sym}, we have $\varphi_l(2a+2)=\varphi_l(2b)$.
Using also $\binom{2l}{2a+2}=\binom{2l}{2b}$, we obtain, for $0\leq b\leq l-1$,
\begin{equation}\label{eq:phi(2b+1)-phi(2b)}
\varphi_l(2b+1)-\varphi_l(2b)=\frac{4^l\binom{2l}{2b}-(4^l-1)\binom{2l}{2b+1}}{2l}.
\end{equation}
On the other hand, \eqref{eq:phi-sum} gives
\[
\sum_{j=0}^{2l}(-1)^j\varphi_l(j)=\varphi_l(0)-\sum_{b=0}^{l-1}(\varphi_l(2b+1)-\varphi_l(2b))=\frac{\varphi_l(0)}{2}.
\]
Thus, by \eqref{eq:phi(2b+1)-phi(2b)}, we can compute
\[
\varphi_l(0)=2\sum_{b=0}^{l-1}(\varphi_l(2b+1)-\varphi_l(2b))=-\frac{4^l}{2l}.
\]
Replacing $b$ by $l-b-1$ in \eqref{eq:phi(2b+1)-phi(2b)} and using the symmetry, we obtain
\[
\varphi_l(2b+1)-\varphi_l(2b+2)=\frac{4^l\binom{2l}{2b+2}-(4^l-1)\binom{2l}{2b+1}}{2l}.
\]
Subtracting this identity from \eqref{eq:phi(2b+1)-phi(2b)}, we get
\[
\varphi_l(2b+2)-\varphi_l(2b)=\frac{4^l}{2l}\left(\binom{2l}{2b}-\binom{2l}{2b+2}\right).
\]
This determines the values of $\varphi_l$ at even indices.
The values at odd indices are then determined by \eqref{eq:phi(2b+1)-phi(2b)}, completing the proof.
\end{proof}
\begin{remark}
Let $p_{2k+1}(x_1,x_2)$ denote Keilthy's polynomial representation of the block-degree-one component of the odd-weight generator $\sigma_{2k+1}$ of the motivic Lie algebra, with the sign conventions and normalization used in \cite[Theorem~5.4]{Keilthy2022}.
Comparing the formula in that theorem with \eqref{eq:explicit_phi}, we obtain
\[
p_{2k+1}(x_1,x_2)=\frac{k}{4^{k-1}}x_1x_2(x_1-x_2)\sum_{j=0}^{2k}\varphi_k(j)x_1^jx_2^{2k-j}.
\]
\end{remark}
\begin{definition}\label{def:eta_rec}
For $\bl=(l_1,\dots,l_r)\in\bbI$ and $(j_1,\dots,j_{r+1})$ with $j_i\in\bbZ_{\geq 2}$, we define $\eta_{\bl}(j_1,\dots,j_{r+1})\in\bbQ$ as follows.
The initial values are given by
\[
\eta_{\varnothing}(j)\coloneqq\begin{cases}
1 & \text{if } j=3,\\
0 & \text{otherwise}.
\end{cases}
\]
For $r\geq 1$, set ${}_{-}\bl\coloneqq (l_2,\dots,l_r)$ and define $\eta_{\bl}$ recursively by
\begin{align*}
&\eta_{\bl}(j_1,\dots,j_{r+1})\\
&=\sum_{\substack{1\leq i\leq r \\ j_i+j_{i+1}\geq 2l_1+3}}(\varphi_{l_1}(j_i-1)-\varphi_{l_1}(j_{i+1}-1))\eta_{{}_{-}\bl}(j_1,\dots,j_{i-1},j_i+j_{i+1}-2l_1-1,j_{i+2},\dots,j_{r+1}).
\end{align*}
\end{definition}
The condition in the sum ensures that every argument of $\eta_{{}_{-}\bl}$ is at least $2$.
By \eqref{eq:explicit_phi} and the recurrence, we have $\eta_{\bl}(j_1,\dots,j_{r+1})\in\bbZ_{(2)}$.
\begin{lemma}\label{lem:eta_nonzero_condition}
We have $\eta_{\bl}(j_1,\dots,j_{r+1})=0$ unless $j_1+\cdots+j_{r+1}=3+2\lvert\bl\rvert+r$.
\end{lemma}
\begin{proof}
At the $i$th step of the recursion, the sum of the block lengths decreases by $2l_i+1$.
A nonzero contribution must eventually reach the initial value $\eta_{\varnothing}(3)$.
\end{proof}
\begin{proposition}\label{prop:eta_T...T}
Let $r$ be a nonnegative integer, let $\bl=(l_1,\dots,l_r)\in\bbI$, and let $v_1,\dots,v_r$ be commuting indeterminates.
Let $w$ be a word satisfying $\block(e_0we_0)=(j_1,\dots,j_{r+1})$ with $j_i\geq2$.
Then
\[
\eta_{\bl}(j_1,\dots,j_{r+1})
=\lambda_{\bl}\cdot[v_1^{2l_1+1}\cdots v_r^{2l_r+1}]
([w]T_{\Phi^-_{\rKZ}(v_1e_0,v_1e_1)}\cdots T_{\Phi^-_{\rKZ}(v_re_0,v_re_1)}(e_1)).
\]
We use the convention $\lambda_{\varnothing}=1$
\end{proposition}
\begin{proof}
We proceed by induction on $r$. When $r=0$, the right-hand side is equal to $[w]e_1$, which is $1$ if $w=e_1$ and $0$ otherwise.
On the other hand, since $\block(e_0we_0)=(j_1)$, the condition $j_1=3$ is equivalent to $w=e_1$.
Hence $\eta_{\varnothing}(j_1)=[w]e_1$.
Let $r\geq 1$.
For each $i$, write $T_{\Phi^-_{\rKZ}(v_ie_0,v_ie_1)}=(T_{\Phi^-_{\rKZ}(v_ie_0,v_ie_1)}-\id)+\id$ and expand the product.
Any term in which $\id$ is chosen for some $i$ has zero coefficient of $v_i^{2l_i+1}$.
Therefore,
\begin{align*}
&[v_1^{2l_1+1}\cdots v_r^{2l_r+1}]
([w]T_{\Phi^-_{\rKZ}(v_1e_0,v_1e_1)}\cdots T_{\Phi^-_{\rKZ}(v_re_0,v_re_1)}(e_1))\\
&=[v_1^{2l_1+1}\cdots v_r^{2l_r+1}]
([w](T_{\Phi^-_{\rKZ}(v_1e_0,v_1e_1)}-\id)\cdots(T_{\Phi^-_{\rKZ}(v_re_0,v_re_1)}-\id)(e_1)).
\end{align*}
Since $\deg_{\rbl}(e_0we_0)=r$, we apply \eqref{eq:T-id_rec} in \cref{lem:T-id} to the right-hand side.
For a decomposition $w=w_0u_1w_1\cdots u_kw_k$, let $w_t$ denote the unique nonempty word among $w_0,\dots,w_k$.
Since
\[
[w_t](\Phi^-_{\rKZ}(v_1e_0,v_1e_1))^{\epsilon(u_t)-\epsilon(u_{t+1})}=v_1^{\lvert w_t\rvert}[w_t](\Phi^-_{\rKZ})^{\epsilon(u_t)-\epsilon(u_{t+1})},
\]
extracting the coefficient of $v_1^{2l_1+1}$ eliminates all terms except those with $\lvert w_t\rvert=2l_1+1$.
For a nonzero term, we have $u_t\neq u_{t+1}$ and $\deg_{\rbl}(u_tw_tu_{t+1})=1$.
Write $e_0we_0=B_1B_2\cdots B_{r+1}$ as a concatenation of alternating blocks, where $\lvert B_i\rvert=j_i$.
We group the terms according to the index $i$ for which $u_tw_tu_{t+1}$ lies in $B_iB_{i+1}$, and express the result as a sum over $i$.
Such an index $i$ must satisfy $j_i+j_{i+1}\geq 2l_1+3$.
By \cref{lem:phi-phi}, we have
\[
\lambda_{l_1}\sum_{\substack{B_iB_{i+1}=pu_tw_tu_{t+1}q\\ u_t,u_{t+1}\in\{e_0,e_1\},\ u_t\neq u_{t+1}\\ \lvert w_t\rvert=2l_1+1\\ \deg_{\rbl}(u_tw_tu_{t+1})=1}}[w_t](\Phi^-_{\rKZ})^{\epsilon(u_t)-\epsilon(u_{t+1})}=\varphi_{l_1}(j_i-1)-\varphi_{l_1}(j_{i+1}-1).
\]
The only pair of equal consecutive letters in $B_iB_{i+1}$ lies inside $u_tw_tu_{t+1}$.
Since $u_t\neq u_{t+1}$, the word $pu_tu_{t+1}q$ is alternating.
Its first letter is the first letter of $B_i$, and its length is $j_i+j_{i+1}-2l_1-1$.
Hence it depends only on $i$, not on the particular decomposition.
Let $w^{(i)}$ denote the resulting word obtained from $w$ by deleting $w_t$.
Then $\block(e_0w^{(i)}e_0)=(j_1,\dots,j_{i-1},j_i+j_{i+1}-2l_1-1,j_{i+2},\dots,j_{r+1})$.
Consequently,
\begin{align*}
&\lambda_{\bl}\cdot[v_1^{2l_1+1}\cdots v_r^{2l_r+1}]([w]T_{\Phi^-_{\rKZ}(v_1e_0,v_1e_1)}\cdots T_{\Phi^-_{\rKZ}(v_re_0,v_re_1)}(e_1))\\
&=\sum_{\substack{1\leq i\leq r\\ j_i+j_{i+1}\geq 2l_1+3}}(\varphi_{l_1}(j_i-1)-\varphi_{l_1}(j_{i+1}-1)) \\
&\qquad\qquad\qquad\cdot\lambda_{{}_{-}\bl}\cdot[v_2^{2l_2+1}\cdots v_r^{2l_r+1}]([w^{(i)}]T_{\Phi^-_{\rKZ}(v_2e_0,v_2e_1)}\cdots T_{\Phi^-_{\rKZ}(v_re_0,v_re_1)}(e_1)).
\end{align*}
Applying the induction hypothesis to each term, the right-hand side becomes the recurrence defining $\eta_{\bl}(j_1,\dots,j_{r+1})$ in \cref{def:eta_rec}.
This proves the identity.
\end{proof}
\begin{corollary}\label{cor:M_eta}
For $N\geq r\geq 1$,
\[
M_{N,r}=(\eta_{\bl}(2k_1+1,\dots, 2k_{r-1}+1, 2k_r+2, 2))_{(\bk,\bl)\in\bbI_{N,r}\times \bbI_{N,r}}.
\]
\end{corollary}
\begin{proof}
The alternating blocks of $e_0\sfw(\sfz(\bk))e_1e_0$ are $B_i=(e_0e_1)^{k_i}e_0$ for each $1\leq i\leq r-1$, $B_r=(e_0e_1)^{k_r+1}$, and $B_{r+1}=e_1e_0$.
Hence,
\begin{equation}\label{eq:block_ewe}
\block(e_0\sfw(\sfz(\bk))e_1e_0)=(2k_1+1,\dots,2k_{r-1}+1,2k_r+2,2).
\end{equation}
Now, by \eqref{eq:adj-Phi=T...T(e_1)} and \cref{prop:eta_T...T}, we obtain
\begin{align*}
\xi_{\bk,\bl}&=\lambda_{\bl}\cdot[v_1^{2l_1+1}\cdots v_r^{2l_r+1}]([\sfw(\sfz(\bk))e_1]\Phi^{-1}_{v_1,\dots,v_r}e_1\Phi_{v_1,\dots,v_r})\\&=\lambda_{\bl}\cdot[v_1^{2l_1+1}\cdots v_r^{2l_r+1}]([\sfw(\sfz(\bk))e_1]T_{\Phi^-_{\rKZ}(v_1e_0,v_1e_1)}\cdots T_{\Phi^-_{\rKZ}(v_re_0,v_re_1)}(e_1))\\&=\eta_{\bl}(2k_1+1,\dots,2k_{r-1}+1,2k_r+2,2),
\end{align*}
which completes the proof.
\end{proof}
\subsection{Proof of \cref{thm:mainG}}
We prove the invertibility of $M_{N,r}$ by reducing a suitable column normalization modulo $2$.
We first establish a congruence for columns indexed by $\bl$ with $l_1=1$.
We then show that the submatrix indexed by indices whose first component is at least $2$ is invertible modulo $2$.
These two results yield the desired invertibility by induction on $r$.
\begin{lemma}\label{lem:odd_vanishing}
Let $r$ be a positive integer, let $\bl=(l_1,\dots,l_r)\in\bbI$, and let $j_1,\dots,j_{r+1}\geq3$ be odd integers.
Then $\eta_{\bl}(j_1,\dots,j_{r+1})\in2^{r+1}\bbZ_{(2)}$.
\end{lemma}
\begin{proof}
By \eqref{eq:explicit_phi}, we have $\varphi_1(2j)\in 2\bbZ_{(2)}$ and $\varphi_l(2j)\in 4\bbZ_{(2)}$ for $l\geq 2$.
We prove the assertion by induction on $r$.
When $r=1$, the recurrence gives
\[
\eta_{(l_1)}(j_1,j_2)=(\varphi_{l_1}(j_1-1)-\varphi_{l_1}(j_2-1))\eta_{\varnothing}(j_1+j_2-2l_1-1)
\]
if $j_1+j_2\geq 2l_1+3$, and $\eta_{(l_1)}(j_1,j_2)=0$ otherwise.
If $l_1\geq 2$, then $j_1-1$ and $j_2-1$ are even by assumption, so both $\varphi_{l_1}(j_1-1)$ and $\varphi_{l_1}(j_2-1)$ belong to $4\bbZ_{(2)}$.
Hence, $\eta_{(l_1)}(j_1,j_2)\in 4\bbZ_{(2)}$.
If $l_1=1$, then \cref{lem:eta_nonzero_condition} implies that $\eta_{(l_1)}(j_1,j_2)\neq 0$ is possible only when $(j_1,j_2)=(3,3)$.
However, $\eta_{(l_1)}(3,3)=0$.

Let $r\geq 2$, and assume that the assertion holds for $r-1$.
It suffices to show that each term on the right-hand side of the recurrence relation for $\eta_{\bl}(j_1,\dots,j_{r+1})$ in \cref{def:eta_rec} belongs to $2^{r+1}\bbZ_{(2)}$. Since $j_i-1$ and $j_{i+1}-1$ are even, we have $\varphi_{l_1}(j_i-1)-\varphi_{l_1}(j_{i+1}-1)\in 2\bbZ_{(2)}$.
Moreover, $j_i+j_{i+1}-2l_1-1$ is also odd, so the induction hypothesis gives
\[
\eta_{{}_{-}\bl}(j_1,\dots,j_{i-1},j_i+j_{i+1}-2l_1-1,j_{i+2},\dots,j_{r+1})\in 2^r\bbZ_{(2)}.
\]
Thus each term in the recurrence relation belongs to $2^{r+1}\bbZ_{(2)}$, completing the proof.
\end{proof}
\begin{lemma}\label{lem:1-equiv_eta}
Let $N\geq r\geq 2$ and $\bk=(k_1,\dots,k_r), \bl=(1,l_2,\dots,l_r)\in\bbI_{N,r}$. Then
\begin{align*}
&\eta_{\bl}(2k_1+1,\dots,2k_{r-1}+1,2k_r+2,2)\\
&\equiv\begin{cases}
-2\eta_{{}_{-}\bl}(2k_2+1,\dots,2k_{r-1}+1,2k_r+2,2) & \text{if } k_1=1,\\
0 & \text{otherwise}
\end{cases}
\pmod{2^r\bbZ_{(2)}}.
\end{align*}
\end{lemma}
\begin{proof}
In fact, we prove the stronger identity
\begin{align*}
&\eta_{\bl}(2k_1+1,\dots,2k_{r-1}+1,2k_r+2,2)\\
&=-2\delta_{k_1,1}\eta_{{}_{-}\bl}(2k_2+1,\dots,2k_{r-1}+1,2k_r+2,2)+3\eta_{{}_{-}\bl}(2k_1+1,\dots,2k_r+1),
\end{align*}
where $\delta$ denotes the Kronecker delta.
Once this identity is established, the assertion follows immediately from \cref{lem:odd_vanishing}.
Set
\[
(j_1,\dots,j_{r+1})=(2k_1+1,\dots,2k_{r-1}+1,2k_r+2,2).
\]
By the recurrence, we have
\[
\eta_{\bl}(j_1,\dots,j_{r+1})=\sum_{i=1}^r(\varphi_1(j_i-1)-\varphi_1(j_{i+1}-1))H_i,
\]
where $H_i\coloneqq\eta_{{}_{-}\bl}(j_1,\dots,j_{i-1},j_i+j_{i+1}-3,j_{i+2},\dots,j_{r+1})$.
Since $2k_r+1\geq3$, we have $\varphi_1(j_r-1)=0$ and $\varphi_1(j_{r+1}-1)=\varphi_1(1)=-3$.
Hence, the term corresponding to $i=r$ is
\[
3H_r=3\eta_{{}_{-}\bl}(2k_1+1,\dots,2k_r+1).
\]
Next, consider the terms with $1\leq i\leq r-1$.
Since $\varphi_1(j_i-1)=-2\delta_{k_i,1}$, we obtain
\begin{align*}
\sum_{i=1}^{r-1}(\varphi_1(j_i-1)-\varphi_1(j_{i+1}-1))H_i&=2\sum_{i=1}^{r-2}(\delta_{k_{i+1},1}-\delta_{k_i,1})H_i-2\delta_{k_{r-1},1}H_{r-1}\\
&=-2\delta_{k_1,1}H_1+2\sum_{i=2}^{r-1}\delta_{k_i,1}(H_{i-1}-H_i).
\end{align*}
Let $2\leq i\leq r-1$. If $k_i=1$, then $j_i=3$, so $H_{i-1}=H_i$. Thus, $\delta_{k_i,1}(H_{i-1}-H_i)=0$ also when $k_i\neq1$.
Consequently,
\[
\sum_{i=1}^{r-1}(\varphi_1(j_i-1)-\varphi_1(j_{i+1}-1))H_i=-2\delta_{k_1,1}H_1.
\]
If $k_1=1$, then $j_1+j_2-3=j_2$, and hence $H_1=\eta_{{}_{-}\bl}(2k_2+1,\dots,2k_{r-1}+1,2k_r+2,2)$.
This completes the proof.
\end{proof}
Set $\bbI'_{N,r}\coloneqq\{\bk=(k_1,k_2,\dots,k_r)\in\bbI_{N,r}:k_1\geq2\}$.
\begin{proposition}\label{prop:M'-invertible}
Let $N\geq r \geq 2$.
Define the matrix $M'_{N,r}$ by
\[
M'_{N,r}\coloneqq(\eta_{\bl}(2k_1+1,\dots,2k_{r-1}+1,2k_r+2,2))_{(\bk,\bl)\in\bbI'_{N,r}\times\bbI'_{N,r}}.
\]
Then the matrix $M'_{N,r}\bmod 2$, obtained by reducing each entry modulo $2\bbZ_{(2)}$, is invertible over $\bbF_2$.
\end{proposition}
\begin{proof}
If $N=r$, then $\bbI'_{N,r}=\varnothing$, and the assertion holds with the convention that the empty matrix is invertible and has determinant $1$.
We therefore assume that $N>r$.
Define $\rho\colon\bbI'_{N,r}\to\bbI'_{N,r}$ by
\[
\rho(l_1,\dots,l_r)=(l_r+1,l_{r-1},\dots,l_2,l_1-1).
\]
Then $\rho$ is an involution.
Let $<_{\rR}$ and $<_{\rL}$ denote the lexicographic orders on $\bbI'_{N,r}$ read from right to left and from left to right, respectively.
Arrange the rows in increasing order with respect to $<_{\rR}$ and the columns in increasing order with respect to $<_{\rL}$.
For example, in $\bbI'_{5,3}$ we have $(3,1,1)<_{\rR}(2,2,1)<_{\rR}(2,1,2)$.
Then $\bk<_{\rR}\bk'$ if and only if $\rho(\bk)<_{\rL}\rho(\bk')$.
With these orderings, we show that $M'_{N,r}$ is upper triangular with all diagonal entries equal to $1$ modulo $2\bbZ_{(2)}$:
\[
M'_{N,r}\equiv
\begin{pmatrix}
1&*&\cdots&*\\
0&1&\cdots&*\\
\vdots&\ddots&\ddots&\vdots\\
0&\cdots&0&1
\end{pmatrix}
\pmod{2\bbZ_{(2)}}.
\]
This will imply that $\det(M'_{N,r})\equiv1\pmod{2\bbZ_{(2)}}$, proving the desired invertibility.
Thus, it suffices to show that
\[
\eta_{\bl}(2k_1+1,\dots,2k_{r-1}+1,2k_r+2,2)\equiv
\begin{cases}
1 & \text{if } \rho(\bl)=\bk,\\
0 & \text{if } \rho(\bl)<_{\rR}\bk
\end{cases}
\pmod{2\bbZ_{(2)}}.
\]
Recall that $\varphi_l(2b)\equiv 0 \pmod{2\bbZ_{(2)}}$ for every integer $b$, that $\varphi_l(2b+1)=0$ whenever $b\geq l$, and that $\varphi_l(2l-1)=\varphi_l(1)=1-4^l\equiv 1 \pmod{2\bbZ_{(2)}}$.

Using the recurrence relation and working modulo $2\bbZ_{(2)}$, we see from $\varphi_l(2b)\equiv 0 \pmod{2\bbZ_{(2)}}$ and \cref{lem:odd_vanishing} that every term except the one obtained by contracting the third and second components from the right is congruent to $0$.
Computing this remaining term gives
\begin{align*}
&\eta_{(l_1,\dots,l_r)}(2k_1+1,\dots,2k_{r-1}+1,2k_r+2,2)\\
&\equiv\begin{cases}
\eta_{(l_2,\dots,l_r)}(2k_1+1,\dots,2k_{r-2}+1,2k_{r-1},2) & \text{if } k_r=l_1-1,\\
0 & \text{if } k_r\geq l_1.
\end{cases}
\end{align*}
Repeating the same calculation, we find that, for $2\leq i\leq r-1$, if
\[
k_r=l_1-1, \quad k_{r-1}=l_2, \quad k_{r-2}=l_3, \quad \dots, \quad k_{r-i+2}=l_{i-1}
\]
then
\begin{align*}
&\eta_{(l_i,\dots,l_r)}(2k_1+1,\dots,2k_{r-i}+1,2k_{r-i+1},2)\\
&\equiv\begin{cases}
\eta_{(l_{i+1},\dots,l_r)}(2k_1+1,\dots,2k_{r-i-1}+1,2k_{r-i},2) & \text{if } k_{r-i+1}=l_i,\\
0 & \text{if } k_{r-i+1}>l_i.
\end{cases}
\end{align*}
Finally, if
\[
k_r=l_1-1, \quad k_{r-1}=l_2, \quad k_{r-2}=l_3, \quad \dots, \quad k_2=l_{r-1}
\]
then the equality $\lvert\bk\rvert=\lvert\bl\rvert=N$ forces $k_1=l_r+1$, that is, $\bk=\rho(\bl)$.
Since
\[
\eta_{(l_r)}(2k_1,2)=\varphi_{l_r}(2k_1-1)-\varphi_{l_r}(1)=4^{l_r}-1\equiv 1 \pmod{2\bbZ_{(2)}},
\]
this completes the verification.
\end{proof}
\begin{proof}[Proof of \cref{thm:mainG}]
The factor $-2$ in \cref{lem:1-equiv_eta} motivates the following column normalization.
For $\bl=(l_1,\dots,l_r)$, define $\kappa(\bl)\coloneqq\max\{i\in\{0,\dots,r-1\}:l_1=\cdots=l_i=1\}$.
Thus $\kappa(\bl)$ is the number of initial entries equal to $1$, capped at $r-1$.
The condition for $i=0$ is understood to be vacuous.
For example, $\kappa(1,2,1)=1$, $\kappa(1,1,2)=2$, and $\kappa(1,1,1)=2$.
In particular, $\kappa(1)=0$.
Define the matrix $\widehat{M}_{N,r}$ over $\bbQ$ by $\widehat{M}_{N,r}\coloneqq(2^{-\kappa(\bl)}\xi_{\bk,\bl})_{(\bk,\bl)\in\bbI_{N,r}\times\bbI_{N,r}}$.
We prove by induction on $r$ that, for every $N\geq r$, all entries of $\widehat{M}_{N,r}$ belong to $\bbZ_{(2)}$ and its reduction modulo $2\bbZ_{(2)}$ is invertible over $\bbF_2$.

First consider the case $r=1$.
We have $\bbI_{N,1}=\{(N)\}$.
Since
\[
\xi_{(N),(N)}=\eta_{(N)}(2N+2,2)=(\varphi_N(2N+1)-\varphi_N(1))\eta_{\varnothing}(3)=4^N-1,
\]
we have $\widehat{M}_{N,1}=M_{N,1}=(4^N-1)$.
Thus, its entry belongs to $\bbZ_{(2)}$, and $\widehat{M}_{N,1}\equiv(1)\pmod{2\bbZ_{(2)}}$.
This establishes the base case.

Let $r\geq2$, and assume that $\widehat{M}_{N',r-1}$ satisfies the assertion for every $N'\geq r-1$.
Fix an arbitrary $N\geq r$, and consider the entries of $\widehat{M}_{N,r}$ indexed by $\bk=(k_1,\dots,k_r),\bl=(l_1,\dots,l_r)\in\bbI_{N,r}$.

For a column with $l_1\geq2$, we have $\kappa(\bl)=0$, so the corresponding columns of $\widehat{M}_{N,r}$ and $M_{N,r}$ coincide.
Their entries $\xi_{\bk,\bl}$ are already known to belong to $\bbZ_{(2)}$.

For an entry with $l_1=k_1=1$, there exists $\varepsilon_{\bk,\bl}\in2^r\bbZ_{(2)}$ such that $\xi_{\bk,\bl}=-2\xi_{{}_{-}\bk,{}_{-}\bl}+\varepsilon_{\bk,\bl}$ by \cref{lem:1-equiv_eta}.
Since $\kappa(\bl)=1+\kappa({}_{-}\bl)$, we obtain $2^{-\kappa(\bl)}\xi_{\bk,\bl}=-2^{-\kappa({}_{-}\bl)}\xi_{{}_{-}\bk,{}_{-}\bl}+2^{-\kappa(\bl)}\varepsilon_{\bk,\bl}$.
By the induction hypothesis, we have $2^{-\kappa({}_{-}\bl)}\xi_{{}_{-}\bk,{}_{-}\bl}\in\bbZ_{(2)}$.
Moreover, $2^{-\kappa(\bl)}\varepsilon_{\bk,\bl}\in2^{r-\kappa(\bl)}\bbZ_{(2)}$, and the definition of $\kappa$ gives $\kappa(\bl)\leq r-1$, so $r-\kappa(\bl)\geq1$.
Consequently, $2^{-\kappa(\bl)}\xi_{\bk,\bl}\in\bbZ_{(2)}$, and
\[
2^{-\kappa(\bl)}\xi_{\bk,\bl}\equiv-2^{-\kappa({}_{-}\bl)}\xi_{{}_{-}\bk,{}_{-}\bl}\pmod{2\bbZ_{(2)}}.
\]

For an entry with $l_1=1$ and $k_1\geq2$, we have $\xi_{\bk,\bl}\in2^r\bbZ_{(2)}$ by \cref{lem:1-equiv_eta}, so $2^{-\kappa(\bl)}\xi_{\bk,\bl}\in2^{r-\kappa(\bl)}\bbZ_{(2)}\subset2\bbZ_{(2)}$.

The map $\bk\mapsto{}_{-}\bk$ gives a bijection between $\{\bk=(k_1,\dots,k_r)\in\bbI_{N,r}:k_1=1\}$ and $\bbI_{N-1,r-1}$.
Thus, grouping both rows and columns according to whether the first component is $1$ or at least $2$, with the former placed first, we obtain
\[
\widehat{M}_{N,r}\equiv
\begin{pmatrix}
-\widehat{M}_{N-1,r-1}&*\\
0&M'_{N,r}
\end{pmatrix}
\pmod{2\bbZ_{(2)}}.
\]
By the induction hypothesis, $\widehat{M}_{N-1,r-1}\bmod{2\bbZ_{(2)}}$ is invertible over $\bbF_2$, and \cref{prop:M'-invertible} shows that $M'_{N,r}\bmod{2\bbZ_{(2)}}$ is also invertible over $\bbF_2$.
Hence, $\widehat{M}_{N,r}\bmod{2\bbZ_{(2)}}$ is invertible over $\bbF_2$.

Consequently, $\det(\widehat{M}_{N,r})\not\equiv0\pmod{2\bbZ_{(2)}}$, so $\det(\widehat{M}_{N,r})\neq0$.
Thus, $\widehat{M}_{N,r}$ is invertible over $\bbQ$.
Since $\widehat{M}_{N,r}$ is obtained from $M_{N,r}$ by multiplying columns by nonzero scalars, $M_{N,r}$ is also invertible over $\bbQ$.
\end{proof}
\section{\texorpdfstring{Proof of \cref{thm:mainF}}{Proof of Theorem \ref*{thm:mainF}}}
We first count the proposed basis elements and then introduce coefficient-extraction maps associated with the family $\Phi_{s;v_1,\dots,v_r}$.
We compute their values on the proposed basis elements and use the resulting vanishing properties, together with the invertibility of $M_{N,r}$, to prove linear independence.
Hirose's dimension bound then completes the proof of \cref{thm:mainF}.
\begin{lemma}\label{lem:gen=d_k}
For every nonnegative integer $k$, we have
\[
\sum_{\substack{a\geq 0,\ N\geq r\geq 1 \\ 2a+2N+r=k}}\sum_{\bk\in\bbI_{N,r}}1+\sum_{\substack{a\geq 0 \\ 2a=k}}1=d_k.
\]
\end{lemma}
\begin{proof}
Let $X$ be an indeterminate. 
We compute the generating function of the left-hand side.
Since
\[
\sum_{N\geq r}\sum_{\bk\in\bbI_{N,r}}X^{2N+r}=\sum_{k_1,\dots,k_r\geq 1}X^{(2k_1+1)+\cdots+(2k_r+1)}=\left(\sum_{m=1}^{\infty}X^{2m+1}\right)^r=\left(\frac{X^3}{1-X^2}\right)^r,
\]
we obtain
\begin{align*}
\sum_{k=0}^{\infty}\left(\sum_{\substack{a\geq 0,\ N\geq r\geq 1 \\ 2a+2N+r=k}}\sum_{\bk\in\bbI_{N,r}}1+\sum_{\substack{a\geq 0 \\ 2a=k}}1\right)X^k&=\sum_{a=0}^{\infty}\sum_{r=1}^{\infty}\sum_{N=r}^{\infty}\sum_{\bk\in\bbI_{N,r}}X^{2a+2N+r}+\sum_{a=0}^{\infty}X^{2a}\\
&=\frac{1}{1-X^2}\left(1+\sum_{r=1}^{\infty}\left(\frac{X^3}{1-X^2}\right)^r\right)\\
&=\frac{1}{1-X^2}\cdot\frac{1}{1-\frac{X^3}{1-X^2}}\\
&=\frac{1}{1-X^2-X^3}.
\end{align*}
The assertion follows by comparing coefficients.
\end{proof}
\begin{definition}
Let $j$ be a nonnegative integer, and let $\bl=(l_1,\dots,l_r)$ be an index.
Define the $\bbQ$-linear map $\Lambda_{j,\bl}\colon\calZ^{\rcf}\to\bbR$ by
\[
\Lambda_{j,\bl}(w\bmod{\calI_{\rCF}})\coloneqq\lambda_{\bl}\cdot[s^{2j}v_1^{2l_1+1}\cdots v_r^{2l_r+1}]\langle w,\Phi_{s;v_1,\dots,v_r}\rangle
\]
for $w\in\calH^0$.
By \cref{lem:pentagon_for_Phi_sv} and \cref{thm:pentagon=confluence}, this map is well-defined.
\end{definition}
\begin{lemma}\label{lem:Lambda_val}
Let $N$ and $r$ be integers satisfying $N\geq r\geq1$, and let $a$ and $j$ be nonnegative integers.
For $\bk\in\bbI_{N,r}$ and an index $\bl=(l_1,\dots,l_d)$, we have
\[
\Lambda_{j,\bl}(\zeta^{\rcf}(2)^a\zeta^{\rcf}_{\calS}(\sfz(\bk)))=\begin{cases}
(-1)^{N+r+1}\zeta(2)^a \xi_{\bk, \bl} & \text{ if } j=a \text{ and } \bl\in\bbI_{N,r}, \\
0 & \text{ if } j=a \text{ and } d>r,\\
0 & \text{ if } j<a.
\end{cases}
\]
Moreover, we have
\[
\Lambda_{j,\bl}(\zeta^{\rcf}(2)^a)=\begin{cases}
\zeta(2)^a & \text{ if } j=a \text{ and } \bl=\varnothing,\\
0 & \text{ otherwise.}
\end{cases}
\]
\end{lemma}
\begin{proof}
We first show that
\begin{equation}\label{eq:zeta(2)s^2}
\langle\zeta^{\rcf}(2),\Phi_{s;v_1,\dots,v_d}\rangle=\zeta(2)s^2.
\end{equation}
If $F=1+(\text{terms of length at least }3)$, then $F^{-1}e_1F=e_1+(\text{terms of length at least }4)$.
In particular, $[w](F\circ G)=[w]G$ for every word $w$ with $\lvert w\rvert\leq2$.
Each series $\Phi^-_{\rKZ}(v_ie_0,v_ie_1)$ has vanishing linear and quadratic terms by its commutator group-like property and \cref{cor:degree2vanish}.
Applying the preceding observation successively to the Ihara factors gives $[e_0e_1]\Phi_{s;v_1,\dots,v_d}=[e_0e_1]\Phi_{\rKZ}(se_0,se_1)=-s^2\zeta(2)$.
On the other hand, $\langle\zeta^{\rcf}(2),\Phi_{s;v_1,\dots,v_d}\rangle=-[e_0e_1]\Phi_{s;v_1,\dots,v_d}$, which proves \eqref{eq:zeta(2)s^2}.

Since $\langle\cdot,\Phi_{s;v_1,\dots,v_d}\rangle$ is an algebra homomorphism, we have
\begin{align*}
\langle\zeta^{\rcf}(2)^a\zeta^{\rcf}_{\calS}(\sfz(\bk)),\Phi_{s;v_1,\dots,v_d}\rangle&=\langle\zeta^{\rcf}(2),\Phi_{s;v_1,\dots,v_d}\rangle^a\langle\zeta^{\rcf}_{\calS}(\sfz(\bk)),\Phi_{s;v_1,\dots,v_d}\rangle\\
&=\zeta(2)^as^{2a}\langle\zeta^{\rcf}_{\calS}(\sfz(\bk)),\Phi_{s;v_1,\dots,v_d}\rangle
\end{align*}
by \eqref{eq:zeta(2)s^2}.
This expression is divisible by $s^{2a}$ and therefore contains no terms of degree less than $2a$ in $s$.
Hence, the coefficient of $s^{2j}$ is zero whenever $j<a$.

When $j=a$, \eqref{eq:sym_cf-Phi} gives
\[
\Lambda_{a,\bl}(\zeta^{\rcf}(2)^a\zeta^{\rcf}_{\calS}(\sfz(\bk)))=(-1)^{N+r+1}\zeta(2)^a\lambda_{\bl}\cdot[v_1^{2l_1+1}\cdots v_d^{2l_d+1}]([\sfw(\sfz(\bk))e_1]\Phi^{-1}_{v_1,\dots,v_d}e_1\Phi_{v_1,\dots,v_d}).
\]
If $\bl\in\bbI_{N,r}$, this is equal to $(-1)^{N+r+1}\zeta(2)^a\xi_{\bk,\bl}$.
By \eqref{eq:adj-Phi=T...T(e_1)} and the argument in the proof of \cref{prop:eta_T...T}, we have
\begin{align*}
&[v_1^{2l_1+1}\cdots v_d^{2l_d+1}]([\sfw(\sfz(\bk))e_1](\Phi^{-1}_{v_1,\dots,v_d}e_1\Phi_{v_1,\dots,v_d}))\\
&=[v_1^{2l_1+1}\cdots v_d^{2l_d+1}]([\sfw(\sfz(\bk))e_1](T_{\Phi^-_{\rKZ}(v_1e_0,v_1e_1)}-\id)\cdots(T_{\Phi^-_{\rKZ}(v_de_0,v_de_1)}-\id)(e_1)).
\end{align*}
By \eqref{eq:block_ewe}, we have $\deg_{\rbl}(e_0\sfw(\sfz(\bk))e_1e_0)=r$.
Since the right-hand side involves $d$ operators of the form $T_{\Phi^-_{\rKZ}(v_ie_0,v_ie_1)}-\id$, when $d>r$ we have $\deg_{\rbl}(e_0\sfw(\sfz(\bk))e_1e_0)<d$ and can apply \eqref{eq:T-id_zero}.
Thus, $\Lambda_{a,\bl}(\zeta^{\rcf}(2)^a\zeta^{\rcf}_{\calS}(\sfz(\bk)))=0$ whenever $d>r$.

By \eqref{eq:zeta(2)s^2}, we have $\langle\zeta^{\rcf}(2)^a,\Phi_{s;v_1,\dots,v_d}\rangle=\zeta(2)^as^{2a}$, and hence
\[
\Lambda_{j,\bl}(\zeta^{\rcf}(2)^a)=\lambda_{\bl}\cdot[s^{2j}v_1^{2l_1+1}\cdots v_d^{2l_d+1}](\zeta(2)^as^{2a}).
\]
The polynomial inside the coefficient extraction does not involve $v_1,\dots,v_d$.
Thus, this expression is zero if $\bl\neq\varnothing$.
If $\bl=\varnothing$, then $\Lambda_{j,\varnothing}(\zeta^{\rcf}(2)^a)=[s^{2j}](\zeta(2)^as^{2a})$, which is zero unless $j=a$, in which case it equals $\zeta(2)^a$.
\end{proof}
\begin{proof}[Proof of \cref{thm:mainF}]
The cases $k=0$ and $k=1$ are immediate, since $\calZ^{\rcf}_0=\bbQ\cdot1$ and $\calZ^{\rcf}_1=0$.
Assume that $k\geq2$.
It suffices to prove that the proposed family is linearly independent.
Indeed, once this is established, \cref{lem:gen=d_k} and \cref{thm:Hirose} imply that it forms a basis.
Suppose, for a contradiction, that there exist rational coefficients $b_{a,\bk}$ and $b_{a,\varnothing}$, not all zero, such that
\begin{equation}\label{eq:linear_comb_contradiction}
\sum_{\substack{a\geq 0,\ N\geq r\geq 1 \\ 2a+2N+r=k}}\sum_{\bk\in\bbI_{N,r}}b_{a,\bk}\zeta^{\rcf}(2)^a\zeta^{\rcf}_{\calS}(\sfz(\bk))+\sum_{\substack{a\geq 0 \\ 2a=k}}b_{a,\varnothing}\zeta^{\rcf}(2)^a=0.
\end{equation}
Let $a_0$ be the smallest value of $a$ for which at least one of the coefficients $b_{a,\bk}$ or $b_{a,\varnothing}$ is nonzero.
Next, set $r_0\coloneqq\max\{\dep(\bk):b_{a_0,\bk}\neq0\}$, where the empty index is included with $\dep(\varnothing)=0$, and set $N_0\coloneqq(k-2a_0-r_0)/2$.

First consider the case $r_0\geq1$.
For an arbitrary $\bl\in\bbI_{N_0,r_0}$, apply $\Lambda_{a_0,\bl}$ to \eqref{eq:linear_comb_contradiction}.
Terms with $a<a_0$ have zero coefficients by the minimality of $a_0$, while terms with $a>a_0$ are annihilated by $\Lambda_{a_0,\bl}$ by \cref{lem:Lambda_val}.
Among the terms with $a=a_0$, those with $r>r_0$ have zero coefficients by the maximality of $r_0$, and those with $r<r_0$ are annihilated because $\dep(\bl)=r_0>r$.
The pure powers of $\zeta^{\rcf}(2)$ are also annihilated, since $\bl\neq\varnothing$.
Thus only the terms with $a=a_0$, $r=r_0$, and $N=N_0$ remain, giving
\[
(-1)^{N_0+r_0+1}\zeta(2)^{a_0}\sum_{\bk\in\bbI_{N_0,r_0}}b_{a_0,\bk}\xi_{\bk,\bl}=0,
\]
that is,
\[
\sum_{\bk\in\bbI_{N_0,r_0}}b_{a_0,\bk}\xi_{\bk,\bl}=0.
\]
This holds for every $\bl\in\bbI_{N_0,r_0}$.
Since $M_{N_0,r_0}$ is invertible by \cref{thm:mainG}, it follows that $b_{a_0,\bk}=0$ for every $\bk\in\bbI_{N_0,r_0}$, contradicting the choice of $a_0$ and $r_0$.

Next, consider the case $r_0=0$.
In this case, $b_{a_0,\varnothing}\neq0$ and $2a_0=k$.
For every term in \eqref{eq:linear_comb_contradiction} with $\bk\neq\varnothing$, we have $a=(k-2N-r)/2<k/2=a_0$, so the minimality of $a_0$ implies that $b_{a,\bk}=0$.
Hence, \eqref{eq:linear_comb_contradiction} reduces to $b_{a_0,\varnothing}\zeta^{\rcf}(2)^{a_0}=0$.
Applying $\Lambda_{a_0,\varnothing}$ and using \cref{lem:Lambda_val}, we obtain $b_{a_0,\varnothing}\zeta(2)^{a_0}=0$.
Since $\zeta(2)>0$, this contradicts $b_{a_0,\varnothing}\neq0$.
This proves the required linear independence and hence \cref{thm:mainF}.
\end{proof}
\section{Proof of main results}
As noted earlier, \cref{thm:mainB} and \cref{thm:mainC} follow immediately from \cref{thm:mainF}.
It remains to prove \cref{thm:mainA}, \cref{thm:mainD}, and \cref{thm:mainE}.
\subsection{Proof of \cref{thm:mainA}}
We transfer the formal expansion \eqref{eq:unique-expansion} to finite multiple zeta values by evaluating at Deligne's $p$-adic associator and reducing modulo $p$.
The key input is a congruence relating symmetric combinations of Deligne's $p$-adic multiple zeta values to finite multiple harmonic sums finally proved by Jarossay.

Let $p$ be a prime.
Let $\Phi^p_{\rKZ}\in\bbQ_p\fps{e_0,e_1}$ and $\Phi^p_{\rDe}\in\bbQ_p\fps{e_0,e_1}$ denote Furusho's $p$-adic KZ associator and Deligne's $p$-adic associator, respectively.
The former was introduced in \cite{Furusho2004} and is related to the latter by
\begin{equation}\label{eq:Deligne_p-adic_associator}
\Phi^p_{\rKZ}=\Phi^p_{\rDe}\circ\Phi^p_{\rKZ}(e_0/p,e_1/p);
\end{equation}
see \cite[Theorem~2.8]{Furusho2007}.
\"{U}nver \cite{Unver2013} proved that $\Phi^p_{\rDe}\in\GRT_1(\bbQ_p)$.

For an index $\bk=(k_1,\dots,k_r)$, we define \emph{Deligne's $p$-adic multiple zeta value} $\zeta_p^{\rDe}(\bk)$ by
\[
\zeta_p^{\rDe}(\bk)\coloneqq(-1)^{\dep(\bk)}[\overleftarrow{\sfw(\bk)}]\Phi^p_{\rDe}.
\]

\begin{theorem}[Chatzistamatiou \cite{Chatzistamatiou2017}, Akagi--Hirose--Yasuda \cite{AkagiHiroseYasuda}, Jarossay \cite{Jarossay}]\label{AHY_congruence}
Let $\bk=(k_1,\dots,k_r)$ be an index. If $p>\lvert\bk\rvert+1$, then all values appearing in the following congruence belong to $\bbZ_p$, and
\[
\sum_{i=0}^r(-1)^{k_{i+1}+\cdots+k_r}\zeta_p^{\rDe}(k_1,\dots,k_i)\zeta_p^{\rDe}(k_r,\dots,k_{i+1})\equiv\zeta^{}_{<p}(\bk)\pmod{p\bbZ_p}.
\]
\end{theorem}
\begin{proof}
The empty-index case is immediate, so assume that $\bk\neq\varnothing$. For $n\geq0$, let $(p\bbZ_p)^{[n]}\subset\bbZ_p$ be the ideal generated by $p^m/m!$ for all $m\geq n$.
By \cite[Theorem~5.3 and Section~5.4]{Chatzistamatiou2017}, we have
\[
[w]\Phi_{\rKZ}^p\in (p\bbZ_p)^{[\lvert w\rvert]}
\]
for every word $w$. An independent proof of this integrality result has also been announced by Akagi, Hirose, and Yasuda \cite{AkagiHiroseYasuda}.
Since $(p\bbZ_p)^{[n]}=p^n\bbZ_p$ when $p>n$, it follows that $[w]\Phi_{\rKZ}^p\in p^{\lvert w\rvert}\bbZ_p$ when $p>\lvert w\rvert$.

Using the inclusions $(p\bbZ_p)^{[u]}(p\bbZ_p)^{[v]}\subset(p\bbZ_p)^{[u+v]}$ and \cref{lem:circ_calc}, we compare coefficients in \eqref{eq:Deligne_p-adic_associator} and argue by induction on word length to obtain
\[
[w]\Phi_{\rDe}^p\in p^{-\lvert w\rvert}(p\bbZ_p)^{[\lvert w\rvert]}.
\]
In particular, $\zeta_p^{\rDe}(\bh)\in\bbZ_p$ for every index $\bh$ with $\lvert\bh\rvert<p$.
Thus, all Deligne's $p$-adic multiple zeta values appearing in the claimed congruence are $p$-adic integers.

The identity conjectured by Akagi, Hirose, and Yasuda and proved by Jarossay \cite{Jarossay} is
\[
\sum_{i=0}^r(-1)^{k_{i+1}+\cdots+k_r}\zeta_p^{\rDe}(k_1,\dots,k_i)\zeta_p^{\rDe,\downarrow}(k_r,\dots,k_{i+1})=\zeta^{}_{<p}(\bk).
\]
In our word convention, the down-arrow values admit the expression
\[
\zeta_p^{\rDe,\downarrow}(\bh)=(-1)^{\dep(\bh)}\sum_{i=0}^{\infty}(-p)^i[\overleftarrow{\sfw(\bh)}e_0^i]\Phi^p_{\rDe}.
\]
Let $\bh$ be nonempty and put $m\coloneqq\lvert\bh\rvert$. For every $i\geq1$, the above coefficient estimate gives
\[
(-p)^i[\overleftarrow{\sfw(\bh)}e_0^i]\Phi^p_{\rDe}\in p^{-m}(p\bbZ_p)^{[m+i]}\subset p^{-m}(p\bbZ_p)^{[m+1]}=p\bbZ_p
\]
whenever $p>m+1$.
Moreover, $(p\bbZ_p)^{[m+i]}$ tends to zero $p$-adically as $i\to\infty$ for $p>2$, so the series converges. Its term with $i=0$ is $\zeta_p^{\rDe}(\bh)$. Hence,
\[
\zeta_p^{\rDe,\downarrow}(\bh)\equiv\zeta_p^{\rDe}(\bh)\pmod{p\bbZ_p}
\]
for $p>\lvert\bh\rvert+1$.
This congruence also holds for the empty index, since both values are $1$.
Applying it to each reversed suffix occurring in the preceding identity and reducing modulo $p\bbZ_p$, we obtain the desired congruence.
\end{proof}
\begin{proof}[Proof of \cref{thm:mainA}]
Since $\Phi^p_{\rDe}\in\GRT_1(\bbQ_p)$, the map $\langle\cdot,\Phi^p_{\rDe}\rangle\colon\calZ^{\rcf}\to\bbQ_p$ is well-defined by \cref{thm:pentagon=confluence}, and we have $\langle\zeta^{\rcf}(2),\Phi^p_{\rDe}\rangle=0$.
Thus, for any nonempty index $\bk$, applying $\langle\cdot,\Phi^p_{\rDe}\rangle$ to the expansion \eqref{eq:unique-expansion} provided by \cref{thm:mainF} gives
\[
\langle\zeta^{\rcf}_{\calS}(\bk),\Phi^p_{\rDe}\rangle=\sum_{\bl\in\bbI^{\rZ}_{\lvert\bk\rvert}}c_{\bk,\bl}\langle\zeta^{\rcf}_{\calS}(\bl),\Phi^p_{\rDe}\rangle.
\]
By \eqref{eq:cf-Phi_general_index}, we have $\langle\zeta^{\rcf}(\bh),\Phi^p_{\rDe}\rangle=\zeta_p^{\rDe}(\bh)$ for every index $\bh$.
Hence \cref{AHY_congruence} applies to the evaluations of the formal symmetric values, and if $p>\lvert\bk\rvert+1$, then
\[
\langle\zeta^{\rcf}_{\calS}(\bk),\Phi^p_{\rDe}\rangle\equiv\zeta^{}_{<p}(\bk)\pmod{p\bbZ_p}.
\]
The same congruence holds for every $\bl\in\bbI^{\rZ}_{\lvert\bk\rvert}$.
Hence, reducing the preceding equality modulo $p$ for all sufficiently large primes $p$ that do not divide the denominator of any coefficient $c_{\bk,\bl}$, and passing to $\calA$, we obtain
\[
\zeta^{}_{\calA}(\bk)=\sum_{\bl\in\bbI^{\rZ}_{\lvert\bk\rvert}}c_{\bk,\bl}\zeta^{}_{\calA}(\bl).
\]
The coefficients $c_{\bk,\bl}$ are rational and independent of $p$, so only finitely many primes are excluded.
This gives the desired expansion with respect to Zhao's indices.
\end{proof}
\subsection{On the Kaneko--Zagier conjecture}
In this subsection, we first prove \cref{thm:mainD}.
We then use the common relation space $\calR_k$ to establish the weightwise Kaneko--Zagier isomorphisms under the corresponding dimension hypotheses.
Combining these isomorphisms with the direct-sum hypotheses yields \cref{thm:mainE}.
We conclude with a remark on compatibility with the depth filtration.
\begin{proof}[Proof of \cref{thm:mainD}]
The proof of \cref{thm:mainA} shows that the coefficients in \eqref{eq:unique-expansion} give the corresponding relations among finite multiple zeta values.
Applying the evaluation map to the same expansion gives the corresponding relations among symmetric multiple zeta values.
Hence every generator of $\calR_k$ in \eqref{eq:def_R_k} belongs to both $\ker Z_{\calA,k}$ and $\ker Z_{\calS,k}$.
The dimension of $\calR_k$ has already been computed in \eqref{eq:dim_of_R_k}.
Moreover, the expansion algorithm described above computes these generators, which form a basis of $\calR_k$.
\end{proof}
\begin{proposition}\label{prop:finiteKanekoZagier_conditional}
Assume that \ref{it:dim_MZV} and \ref{it:dim_FMZV} hold.
Then \cref{conj:KanekoZagier_weight} is true.
\end{proposition}
\begin{proof}
Fix $k\geq2$.
The case $k=2$ is immediate, since both spaces are zero.
Assume that $k\geq3$.
Combining the known upper bounds with the assumed lower bounds, we obtain $\dim_{\bbQ}\calZ_{\calA,k}=d_{k-3}$ and $\dim_{\bbQ}(\calZ_k/\zeta(2)\calZ_{k-2})=d_k-d_{k-2}=d_{k-3}$.
It follows that $\dim_{\bbQ}\ker Z_{\calA,k}=\dim_{\bbQ}\ker Z_{\calS,k}=2^{k-1}-d_{k-3}$, and hence \cref{thm:mainD} gives $\ker Z_{\calA,k}=\calR_k=\ker Z_{\calS,k}$.

Let $\overline{Z_{\calA,k}}$ and $\overline{Z_{\calS,k}}$ denote the isomorphisms induced by $Z_{\calA,k}$ and $Z_{\calS,k}$ on the common quotient $\calH^1_k/\calR_k$.
Then
\[
\varphi_{\rKZ,k}\coloneqq\overline{Z_{\calS,k}}\circ(\overline{Z_{\calA,k}})^{-1}\colon\calZ_{\calA,k}\overset{\sim}{\longrightarrow}\calZ_k/\zeta(2)\calZ_{k-2}
\]
is the desired isomorphism.
By construction, this isomorphism sends $\zeta^{}_{\calA}(\bk)$ to $\zeta^{\sh}_{\calS}(\bk)\bmod{\zeta(2)\calZ_{k-2}}$ for every $\bk\in\bbI_k$.
\end{proof}
\begin{lemma}\label{lem:direct_sum_modz2}
Assume that \ref{it:direct_sum_MZV} holds.
Then the natural surjective map
\[
\bigoplus_{k=0}^{\infty}(\calZ_k/\zeta(2)\calZ_{k-2})\twoheadrightarrow\calZ/\zeta(2)\calZ
\]
given by summing the components is an isomorphism of graded $\bbQ$-algebras.
\end{lemma}
\begin{proof}
The map is a surjective homomorphism of graded $\bbQ$-algebras.
It remains to prove injectivity.
Suppose that a finite sum
\[
\sum_k(a_k\bmod{\zeta(2)\calZ_{k-2}})
\]
lies in the kernel of this map, where $a_k\in\calZ_k$.
Then there exists a finite sum $\sum_j b_j$ with $b_j\in\calZ_j$ such that
\[
\sum_k a_k=\zeta(2)\sum_j b_j.
\]
Set $b_j=0$ for all indices not occurring in the sum, including $j<0$.
By \ref{it:direct_sum_MZV}, comparing the components of weight $k$ gives $a_k=\zeta(2)b_{k-2}\in\zeta(2)\calZ_{k-2}$ for every $k$.
Thus the original element of the kernel is zero, proving injectivity.
\end{proof}
\begin{proof}[Proof of \cref{thm:mainE}]
In weight $0$, let $\varphi_{\rKZ,0}$ be the identity map on $\bbQ$.
In weight $1$, let $\varphi_{\rKZ,1}$ be the unique map between the zero spaces $\calZ_{\calA,1}$ and $\calZ_1/\zeta(2)\calZ_{-1}$.
First, finite multiple zeta values and symmetric multiple zeta values both satisfy the stuffle product formula, which preserves weight; see \cite[Proposition~4 in Section~6 and Proposition~8 in Section~8]{KanekoZagier}.
Consequently, assuming \ref{it:dim_MZV} and \ref{it:dim_FMZV}, the weightwise isomorphisms in the proof of \cref{prop:finiteKanekoZagier_conditional} combine to give an isomorphism of graded $\bbQ$-algebras
\[
\bigoplus_{k\geq 0}\varphi_{\rKZ,k}\colon\bigoplus_{k\geq 0}\calZ_{\calA,k}\overset{\sim}{\longrightarrow}\bigoplus_{k\geq 0}(\calZ_k/\zeta(2)\calZ_{k-2}).
\]
Since each $\varphi_{\rKZ,k}$ preserves the index, the stuffle product formulas imply that their direct sum preserves products.
Thus, if we further assume \ref{it:direct_sum_MZV} and \ref{it:direct_sum_FMZV}, composing with the isomorphism in \cref{lem:direct_sum_modz2} yields the desired isomorphism of graded algebras
\[
\varphi_{\rKZ}\colon\calZ_{\calA}\overset{\sim}{\longrightarrow}\calZ/\zeta(2)\calZ.\qedhere
\]
\end{proof}
\begin{remark}\label{rem:depth_KZconj}
For integers $k>0$ and $d\geq 0$, we set
\begin{align*}
\calD_d\calZ_k&\coloneqq\Span_{\bbQ}\{\zeta(\bk) : \bk\in\bbI_k \text{ is admissible, } \dep(\bk)\leq d\},\\
\calD_d(\calZ_k/\zeta(2)\calZ_{k-2})&\coloneqq(\calD_d\calZ_k+\zeta(2)\calZ_{k-2})/\zeta(2)\calZ_{k-2},\\
\calD_d\calZ_{\calA,k}&\coloneqq\Span_{\bbQ}\{\zeta^{}_{\calA}(\bk) : \bk\in\bbI_k,\ \dep(\bk)\leq d+1\},\\
\calS_{k,d}&\coloneqq\Span_{\bbQ}\{\zeta^{}_{\calS}(\bk) : \bk\in\bbI_k,\ \dep(\bk)\leq d+1\}.
\end{align*}
In \cref{conj:KanekoZagier}, we interpret compatibility of $\varphi_{\rKZ}$ with the depth filtration as the condition that
\[
\varphi_{\rKZ,k}(\calD_d\calZ_{\calA,k})=\calD_d(\calZ_k/\zeta(2)\calZ_{k-2})
\]
for all $k$ and $d$. Under the assumptions in \cref{thm:mainE}, $\varphi_{\rKZ}$ is an isomorphism, so it remains to prove that
\[
\calS_{k,d}=\calD_d(\calZ_k/\zeta(2)\calZ_{k-2}).
\]
The inclusion $\calS_{k,d}\subset\calD_d(\calZ_k/\zeta(2)\calZ_{k-2})$ was established in \cite[Theorem~4 in Section~8]{KanekoZagier}.
Thus, establishing the reverse inclusion $\calS_{k,d}\supset\calD_d(\calZ_k/\zeta(2)\calZ_{k-2})$ would strengthen the conclusion of \cref{thm:mainE} to the full Kaneko--Zagier conjecture.

In the proof of the corollary immediately following \cite[Theorem~4 in Section~8]{KanekoZagier}, Kaneko and Zagier report a depth-refined version of Yasuda's theorem \cref{thm:Yasuda}, attributed to a personal communication from Yasuda in 2016.
Such a refinement, in the filtered form above, would give the reverse inclusion and hence compatibility with the depth filtration.
A proof of this refinement is not given there, and we do not use it in the present paper.
\end{remark}
\subsection*{Acknowledgement}
The author would like to express his sincere gratitude to Professor Hidekazu Furusho for kindly answering questions about his work.
\subsection*{Use of AI}
ChatGPT (OpenAI, GPT-6 Astra Pro) was used for mathematical discussions and English-language editing.
The final manuscript and all mathematical arguments were written and independently verified by the author, who takes full responsibility for the content.

\end{document}